\documentclass[11pt,reqno]{amsart}

\usepackage{amsmath,mathtools}
\usepackage{graphicx}
\usepackage{amsfonts}
\usepackage{amssymb}
\usepackage{setspace}
\usepackage{color}
\usepackage[font=small]{caption}
\usepackage[font=footnotesize]{subcaption} 
\usepackage[top=1in, bottom=1in, left=1.2in, right=1.2in]{geometry}
\usepackage{afterpage}
\usepackage{bm}
\usepackage{multicol}
\usepackage{multirow}
\usepackage{isomath}

\usepackage{float}
\usepackage{tikz}
\usetikzlibrary{patterns}
\usetikzlibrary{math}
\usepackage{algorithm}
\usepackage{algorithmicx}
\usepackage{algpseudocode}
\usepackage{amsthm}

\tikzset{global scale/.style={
    scale=#1,
    every node/.append style={scale=#1}
} }

\newtheorem{theorem}{Theorem}[section]

\newtheorem{remark}[theorem]{Remark}

\def \real{\mathbb{R}}        
\def \complex{\mathbb{C}}     
\def \1:#1{1,\dots,#1}

\def \eps{\varepsilon}           
\def \O(#1){\mathcal{O}\left(#1\right)}  
\def \norm(#1){\left\| #1 \right\|} 
\def \abs(#1){\left| #1 \right|}  

\def \fig{Fig.~}  
 
\def \algo{Algorithm~} 
\def \equ(#1){equation~(#1)}
\def \inequ(#1){inequality~(#1)}

\graphicspath{{fig}}

\title[Recursive sparse LU decomposition]{
	Recursive sparse LU decomposition based on nested dissection and low rank approximations}

\author{Xuanru Zhu}
\address{School of Mathematical Sciences, Zhejiang University,
	Hangzhou, Zhejiang 310027, China}
\email{zhuxuanru@126.com}

\author{Jun Lai}
\address{School of Mathematical Sciences, Zhejiang University,
	Hangzhou, Zhejiang 310027, China}
\email{laijun6@zju.edu.cn}

\subjclass[2020]{65F05, 65F55, 65N30, 35J05}

\keywords{Hierarchical matrices; nested dissection; low rank approximations; fast direct solver}

\date{\today}

\begin{document}

\begin{abstract}
When solving partial differential equations (PDEs) using finite difference or finite element methods, efficient solvers are required for handling large sparse linear systems. In this paper, a recursive sparse LU decomposition for matrices arising from the discretization of linear PDEs is proposed based on the nested dissection and low rank approximations.  The matrix is reorganized based on the nested structure of the associated graph. After eliminating the interior vertices at the finest level, dense blocks on the separators are hierarchically sparsified using low rank approximations.  To efficiently skeletonize these dense blocks, we split the separators into segments and introduce a hybrid algorithm to extract the low rank structures based on a randomized algorithm and the fast multipole method.  The resulting decomposition yields a fast direct solver for sparse matrices, applicable to both symmetric and non-symmetric cases.  Under a mild assumption on the compression rate of dense blocks, we prove an $\O(N)$  complexity for the fast direct solver. Several numerical experiments are provided to verify the effectiveness of the proposed method. 
\end{abstract}
\maketitle

\section{Introduction}
    Consider the following partial differential equations:
\begin{equation}\label{generalpde}
	\begin{cases}
		\mathcal{L}u = f \mbox{ in } \Omega,\\
		\mathcal{T}u = g\mbox{ on }\partial \Omega,
	\end{cases}
\end{equation}
where $\mathcal{L}$ is a linear differential operator, $\mathcal{T}$ is a boundary operator, $f$ and $g$ are the source and boundary functions, respectively. The computational domain $\Omega\in \mathbb{R}^d$ is typically assumed to be simply connected with a Lipschitz continuous boundary. Such kind of PDEs,  examples including Laplace and Helmholtz equations with Dirichlet or Neumann boundary conditions \cite{cave2014, Lai2015194}, commonly arise  in many important applications. When equation \eqref{generalpde} is solved by finite difference (FD) or finite element method (FEM), it results in a linear system: 
\begin{equation}\label{equ-main}
\mathcal{A} x=b,
\end{equation}
where matrix $\mathcal{A}$ is large but sparse. For many applications  involving optimizations with PDE constraints, such as optimal design \cite{MaLin2024} and inverse problems  \cite{inverseP2022}, a key challenge is how to solve equation \eqref{equ-main} efficiently.

One common approach for solving equation \eqref{equ-main} involves iterative methods based on Krylov subspaces, such as  Conjugate Gradient (CG) \cite{MatrixComputation} or Generalized Minimal Residual (GMRES) \cite{GMRES}, which typically costs near $\O(N)$ complexity for the matrix-vector product and can complete the iterations in $\O(1)$ steps when paired with an effective preconditioner. However, the construction of such a preconditioner is often highly nontrivial and many of the existing techniques are problem-specific. A common general-purpose preconditioner  is the incomplete LU (ILU) \cite{ILU}, which begins with the classical LU factorization and ignores some of the \textit{fill-in}s below a given threshold. Despite its utility, ILU itself can be computationally expensive, since the number of fill-ins may significantly exceed that of the original sparse matrix. Another well-known approach is the \textit{geometric} multigrid method \cite{Multigrid}, which is an iterative algorithm based on a sequence of mesh discretizations at multiple scales. It makes use of the smoothing effect of the classical Jacobi or Gauss-Seidel iterations and takes the solution from the coarse grid as the initial guess for the finer grid, leading to rapid convergence for non-oscillatory solutions. However, for equations with oscillatory solutions, such as the Helmholtz and Maxwell equations \cite{cave2014}, the convergence rate of multigrid methods tends to deteriorate significantly. In addition, while iterative methods can be very effective for problems with a single right hand side, they become inefficient when multiple right hand sides need to be solved, in which case direct solvers based on the factorization of matrices are highly desirable.

Given the mesh information of PDEs, nested dissection (ND) has emerged as a natural approach for factorizing the matrix from FD or FEM discretization \cite{ND-intro_first}. The basic idea of ND is to recursively partition the graph $G$ associated with the sparse matrix $\mathcal{A}$ into two (or more) disjoint subgraphs using a series of separators, continuing this process  until each subgraph at the leaf level contains $\mathcal{O}(1)$ vertices. The method then involves factorizing matrix blocks in a bottom-up manner, starting from the finest level.   It is a typical \textit{divide-and-conquer} approach that effectively makes use of the local connectivity properties of the vertices in the FD or FEM graph. However, as the elimination progresses, the system becomes increasingly dense  due to the full connectivity of the vertices along the separators. This typically leads to a cost of $\mathcal{O}(N^{3/2})$ by a straightforward LU decomposition for separators of size $\mathcal{O}(N^{1/2})$ in two dimensions and  $\mathcal{O}(N^{2})$ for separators of size $\mathcal{O}(N^{2/3})$ in three dimensions \cite{ND-timecost}. 

To enhance the efficiency, methods based on low rank approximations have been introduced to accelerate the computation of dense blocks~\cite{Intro-solveEq}, utilizing the structure of \textit{Hierarchical} matrix($\mathcal{H}$-matrix) \cite{Hierarchical}. These approaches make use of the fact that the underlying Green's functions for the PDEs, while often oscillatory, are smoothly decaying in their amplitudes, leading to far-field interactions that are low rank up to a small tolerance. Examples of $\mathcal{H}$-matrices include hierarchically semiseparable (HSS) matrices \cite{HSS}, hierarchical off-diagonal low rank matrices(HODLR) \cite{HODLR}, $H^2$ matrices \cite{H^2} etc.  By hierarchically extracting low rank structures and factorizing them recursively, it could lead to an approximate inverse factorization with asymptotically linear complexity. Due to these advantages, they have been widely used in solving the dense linear systems derived from the discretizations of boundary integral equations \cite{FMM-lecture}. However, when it comes to linear systems from FD or FEM discretizations, the lack of exploitation of unique data structure existed in FD or FEM graphs results in significant overhead in time complexity. 

Recently, Cambier et al.\cite{spaND} proposed a fast factorization algorithm called \textit{spaND}, which stands for \textit{sparsified Nested Dissection}, for the sparse matrices based on nested dissection. This algorithm employs low rank approximations to directly sparsify the separators in the graph without introducing additional fill-ins. Unlike other fast methods for ND that store large blocks in a low rank format, and then compress these dense blocks via fast $\mathcal{H}$-algebra, spaND sparsifies and eliminates the separators right from the beginning. The algorithm is particularly designed for symmetric positive definite (SPD) matrices and yields a hierarchically sparse $GG^\top$ decomposition with nearly linear complexity. Such decomposition is then used as a preconditioner to construct an efficient iterative solver. SpaND fully explores the unique  structures in FD or FEM graphs by locally sparsifying dense blocks on the separators. However, it still has certain limitations: it is restricted to symmetric matrices, and the efficiency of low rank structure extraction via QR decomposition is relatively low. Meanwhile, although it is a common strategy to use the sparse factorization as a preconditioner in iterative solvers, it does not fully take the advantages of sparse representation, as compared to the construction of fast direct solver, in the case when multiple right hand sides need to be solved. 

Our work builds on the framework of spaND with LU decomposition to develop a method called spaLU, which stands for sparse LU decomposition. It begins with the nested dissection of a finite element graph, and recursively sparsifies the dense blocks on the separators through low rank approximations, ultimately resulting in a recursive sparse LU decomposition. Compared to spaND, spaLU offers the following enhancements:

\begin{enumerate}
\item  Based on the sparse LU decomposition, we extend the algorithm from SPD to any invertible sparse matrix $\mathcal{A}$ derived from FD or FEM discretization without significantly increasing the computational cost.
  
\item We propse a hybrid sampling method from randomized sampling and FMM that accelerates the extraction of low rank structures, providing stable sparsification of dense blocks with a controllable error.

\item We present a fast direct solver by utilizing the recursive LU decomposition and demonstrate its linear complexity under a mild assumption on the sparsification rate. It is particularly effective for solving equations with multiple right hand sides.
\end{enumerate}	

The paper is organized as follows: Section 2 introduces the notations for nested dissection and provides an overview of recursive sparse LU decomposition. Section 3 details the construction of data structure in nested dissection, focusing on separators and segments. Section 4 proposes an efficient sparsification of segments based on a hybrid interpolative decomposition algorithm, which combines randomized sampling and the field separation idea from FMM. Section 5 presents the complexity analysis of the algorithm. Numerical experiments are provided in Section 6, and the paper is concluded in Section 7.
\section{Preliminaries}

In this section, we introduce the notations and outline the main idea of recursive sparse LU decomposition method. For clarity, we focus on the two dimensional FEM graph, with the extension to three dimensional PDEs to be addressed in future work.

\subsection{Graph and elimination}
\label{section-graph and elimination}

Graph structure plays an important role in the computation of sparse matrices \cite{ND-graph_for_eliminate}. Direct methods such as LU, Cholesky, and QR factorizations can be conveniently illustrated by graphs during the elimination process.

Specifically, a large sparse matrix $\mathcal{A}$ of order $n$ corresponds to an undirected graph $G=(V,E)$ consisted of $n$ \textit{vertices} $V = \{v_i, \ i=1,\dots,n\}$ and \textit{edges} $E=\{e_{ij}\}$. An edge $e_{ij}$ connects vertex $v_i$ to $v_j$ if and only if $\mathcal{A}_{ij}$ is nonzero, making vertex $v_j$ the adjacent vertex of $v_i$.
Denote $\mathrm{Adj}(v_i)$ the set of all adjacent vertices of $v_i$, and $\mathrm{Deg}(v_j)$ the \textit{degree} of a vertex $v_i$, which is the number of its adjacent vertices. The adjacent vertices of a vertex set $V$ are denoted by $\mathrm{Adj}(V) = \bigcup_{v_i \in V} \mathrm{Adj}(v_i)$.
Generally, the graph associated with a matrix $\mathcal{A}$ is constructed based on the non-zero elements of $\mathcal{A}$ \cite{ND-graph_for_eliminate}. However, for matrix obtained from the discretization of PDEs, the corresponding graph can be constructed easily from the discretization mesh. In particular, the matrix graph for FEM with linear or bilinear elements is identical to the discretization mesh, often referred to as the finite element graph.
		
When eliminating vertex $v$ from the graph $G$ through LU factorization, new non-zero elements, known as \textit{fill-ins}, are generally introduced into the updated matrix. These fill-ins correspond to new edges in the graph. Define the deficiency $\mathrm{Def}(v)$ as the set of new edges after eliminating $v$:
     \begin{equation}
     	\mathrm{Def}(v) = \{e_{ij} \, | \, v_i, v_j \in \mathrm{Adj}(v), \mbox{ and } e_{ij} \notin E \}.
     \end{equation}
     The new graph after eliminating $v$ (also referred to as $v$-\textit{elimination}) is given by
	\begin{equation}
		G_v = (V_v, E_v)\mbox{ with }V_v = V \backslash \{v\}, \  
		 E_v=\{e_{ij} \in E| v_i \mbox{ or } v_j\ne v \} \cup \mathrm{Def}(v).
	\end{equation}
 
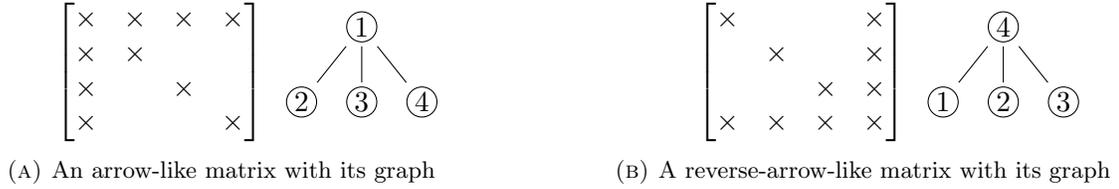
\begin{figure}[tb] 
	\centering
	\begin{subfigure} [t] {0.48\textwidth}
		\centering
		\begin{minipage}[c]{0.22\textwidth}
			\centering
			\begin{equation*}
				\begin{bmatrix}
					\times & \times & \times & \times \\
					\times & \times & & \\
					\times & & \times & \\
					\times & & & \times
				\end{bmatrix}
			\end{equation*}
		\end{minipage}
		\quad\quad\quad
		\begin{minipage}[c]{0.17\textwidth}
			\centering
			\begin{tikzpicture}
				\tikzmath{ \len = 0.8; };
				\node (v1) at (\len,1) {1};
				\draw (v1) circle [radius=0.2];
				
				\node (v2) at (0,0) {2};
				\node (v3) at (\len,0) {3};
				\node (v4) at (\len*2,0) {4};
				\draw (v2) circle [radius=0.2];
				\draw (v3) circle [radius=0.2];
				\draw (v4) circle [radius=0.2];
				
				\draw (v1)--(v2);
				\draw (v1)--(v3);
				\draw (v1)--(v4);
			\end{tikzpicture}
		\end{minipage}
		
		\caption{An arrow-like matrix with its graph}
		\label{fig-arrow-like 1}
	\end{subfigure}
	\quad
	\begin{subfigure} [t] {0.48\textwidth}
		\centering
		\begin{minipage}[c]{0.22\textwidth}
			\centering
			\begin{equation*}
				\begin{bmatrix}
					\times & & & \times \\
					& \times & & \times \\
					& & \times & \times \\
					\times & \times & \times & \times
				\end{bmatrix}
			\end{equation*}
		\end{minipage}
		\quad\quad\quad
		\begin{minipage}[c]{0.17\textwidth}
			\centering
			\begin{tikzpicture}
				\tikzmath{ \len = 0.8; };
				\node (v1) at (\len,1) {4};
				\draw (v1) circle [radius=0.2];
				
				\node (v2) at (0,0) {1};
				\node (v3) at (\len,0) {2};
				\node (v4) at (\len*2,0) {3};
				\draw (v2) circle [radius=0.2];
				\draw (v3) circle [radius=0.2];
				\draw (v4) circle [radius=0.2];
				
				\draw (v1)--(v2);
				\draw (v1)--(v3);
				\draw (v1)--(v4);
			\end{tikzpicture}
		\end{minipage}
		
		\caption{A reverse-arrow-like matrix with its graph}
		\label{fig-arrow-like 2}
	\end{subfigure}
	
	\caption{ A reordering of matrix elements. (a): In the arrow-like matrix, a large number of fill-ins will be introduced in the  LU elimination. (b): In the reordered reverse-arrow-like matrix, the LU elimination  creates a minimal number of fill-ins.}
	\label{figur1}
\end{figure}

 It is evident that different eliminations orders cause fill-ins of different sizes. For example, as shown in \fig \ref{figur1}, eliminating an arrow-like matrix in \fig \ref{fig-arrow-like 1} in ascending order (from row $1$ to row $4$) by LU decomposition results in a large number of fill-ins in the matrix factors $L$ and $U$. In contrast, eliminating the same matrix in the reversed order, as shown in \fig \ref{fig-arrow-like 2}, results in a minimal number of fill-ins. Therefore, it is crucial to select an elimination order such that the degree of $\mathrm{Def}(v)$ is minimized after the  $v$-\textit{elimination}. However, finding an optimal elimination order of a given matrix is usually NP-complete \cite{ND-NPC}. One heuristic approach is the \textit{Minimum Degree Ordering} \cite{MinDegreeOrder}, which rearranges the elements in ascending order of their degrees, but in most cases it is still not optimal. Nested dissection offers an approach to address this issue by reordering the vertices to confine the number of fill-ins introduced by elimination to well-separated subgraphs. This approach localizes newly generated edges, with interactions spreading gradually to the entire graph.
 
\subsection{Nested dissection}  
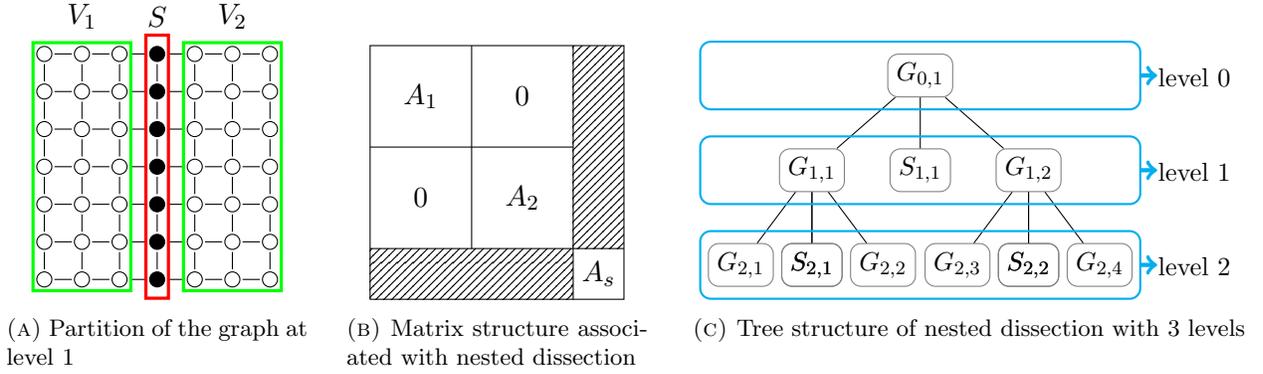
\begin{figure}[tb] 
    \centering
    \begin{subfigure} [t] {0.24\textwidth}
        \centering
        \begin{tikzpicture}
            \tikzmath{ \len = 0.5; };
            \foreach \j in {1,2,...,7}
            {
                \foreach \i in {1,2,...,7}
                {
                    \node (a\j\i) at (\len*\i-\len,\len*\j-\len) {};
                    \draw (a\j\i) circle [radius=0.1];
                }
                \foreach \i/\t in {2/1,3/2,4/3,5/4,6/5,7/6}
                {
                    \draw (a\j\t)--(a\j\i);
                }
            }
            \foreach \j in {1,2,...,7}
            {
                \foreach \i/\t in {2/1,3/2,4/3,5/4,6/5,7/6}
                {
                    \draw (a\t\j)--(a\i\j);
                }
            }
            \foreach \i in {1,2,...,7}
            {
                \fill[black] (a\i 4) circle [radius=0.1];
            }
            
            \draw[green][very thick] (-0.15,-0.15) rectangle (\len*2+0.15,\len*6+0.15);
            \draw [green][very thick] (\len*4-0.15,-0.15) rectangle (\len*6+0.15,\len*6+0.15);
            \draw [red][very thick] (\len*3-0.15,-0.25) rectangle (\len*3+0.15,\len*6+0.25);
            
            \node at (\len*1,\len*6+0.5) {$V_1$};
            \node at (\len*3,\len*6+0.5) {$S$};
            \node at (\len*5,\len*6+0.5) {$V_2$};
        \end{tikzpicture}

        \caption{Partition of the graph at level 1}
        \label{fig-nested graph}
    \end{subfigure}
    \quad
    \begin{subfigure} [t] {0.24\textwidth}
        \centering
        \begin{tikzpicture} [scale=0.9,
            pattern_slash/.style={pattern=north east lines}]
            \draw (3.75,3.75) rectangle (0,0);
                
            \draw (3,0) -- (3,3.75);
            \draw (0,0.75) -- (3.75,0.75);
            
            \draw (1.5,0.75) -- (1.5,3.75);
            \draw (0,2.25) -- (3,2.25);
            
            \fill[pattern_slash] (0,0) rectangle (3,0.75);
            \fill[pattern_slash] (3.75,3.75) rectangle (3,0.75);
            
            \node at (0.75,3) {$A_1$};
            \node at (2.25,1.5) {$A_2$};
            \node at (3.375,0.375) {$A_s$};
            \node at (2.25,3) {0};
            \node at (0.75,1.5) {0};
        \end{tikzpicture}
        
        \caption{Matrix structure associated with nested dissection}
        \label{fig-nested matrix}
    \end{subfigure}
    \quad
    \begin{subfigure} [t] {0.45\textwidth}
        \centering
        \begin{tikzpicture} [global scale = 0.9]
            \node [rectangle,rounded corners,draw=gray] (a11) at (0,0) {$G_{0,1}$};
            \foreach \i in {1,2}
            {
                \node [rectangle, rounded corners, draw=gray] (a2\i) at (-4.8+3.2*\i,-1.4) {$G_{1,\i}$};
                \draw (a11) -- (a2\i);
            }
            \foreach \i in {1,2}
            {
                \tikzmath{ int \I; \I = \i; };
                \node [rectangle, rounded corners, draw=gray] (a3\I) at (-4.75+2.1*\I,-2.8) {$G_{2,\I}$};
                \draw (a21) -- (a3\I);
                \node [rectangle, rounded corners, draw=gray] (a3c) at (-1.6,-2.8) {$S_{2,1}$};
                \draw (a21) -- (a3c);
                \tikzmath{ int \I; \I = 2 + \i; };
                \node [rectangle, rounded corners, draw=gray] (a3\I) at (-5.75+2.1*\I,-2.8) {$G_{2,\I}$};
                \draw (a22) -- (a3\I);
                \node [rectangle, rounded corners, draw=gray] (a3d) at (1.6,-2.8) {$S_{2,2}$};
                \draw (a22) -- (a3d);
            }
            \node [rectangle, rounded corners, draw=gray] (a2c) at (0,-1.4) {$S_{1,1}$};
            \draw (a11) -- (a2c);
            
            \draw[thick, rounded corners, cyan] (-3.25,0.5) rectangle (3.25,-0.5);
            \draw[->] [ultra thick, cyan] (3.25,0)--(3.5,0);
            \node at (4.05,-0.025) {level 0};
            \draw[thick, rounded corners, cyan] (-3.25,-0.9) rectangle (3.25,-1.9);
            \draw[->] [ultra thick, cyan] (3.25,-1.4)--(3.5,-1.4);
            \node at (4.05,-1.425) {level 1};
            \draw[thick, rounded corners, cyan] (-3.25,-2.3) rectangle (3.25,-3.3);
            \draw[->] [ultra thick, cyan] (3.25,-2.8)--(3.5,-2.8);
            \node at (4.05,-2.825) {level 2};
        \end{tikzpicture}
        
        \caption{Tree structure of nested dissection with 3 levels}
        \label{fig-split tree}
    \end{subfigure}
    
    \caption{An illustration of nested dissection. }
\end{figure}

Given the graph $G$ of a matrix $\mathcal{A}$, the nested dissection algorithm \cite{ND-intro_lite, ND-intro_first} partitions the graph $G = (V,E)$ into three parts: disjoint subgraphs $G_1 = (V_1, E_1)$, $G_2 = (V_2, E_2)$, and the \textit{separator} $S$, such that
\begin{gather}
    \begin{cases}
        V = V_1 \cup V_2 \cup S,  \\
        E_1, E_2 \subset E, \\
        \{e_{ij} \in E \, | \, v_i \in V_1 , \, v_j \in V_2 \} = \emptyset,
        \end{cases}
\end{gather}
as illustrated in \fig \ref{fig-nested graph}. The algorithm then reorders the rows and columns of $\mathcal{A}$ so that elements correspond to $S$ are placed in the end of matrix $\mathcal{A}$, following the blocks associated with $G_1$ and $G_2$, as shown in \fig \ref{fig-nested matrix}. The LU elimination of the two blocks associated with $G_1$ and $G_2$ does not affect each other. In addition, subgraphs $G_1$ and $G_2$  can be recursively divided in a nested manner until there are $\mathcal{O}(1)$ vertices in each subgraph at the leaf level. This procedure can be represented by a splitting tree, as illustrated in \fig \ref{fig-split tree}. 

The factorization procedure for matrix $\mathcal{A}$ begins by eliminating the elements corresponding to interior vertices in each subgraphs at the leaf level, and then successively eliminating the separators from the finest level to the coarsest.  Denote $A$ as the reordered matrix $\mathcal{A}$ with nested order $I$. The nested dissection algorithm converts equation \eqref{equ-main} into 
\begin{equation}
    A y = \tilde{b}
    \label{equ-Ax=b}
\end{equation}
with
\begin{equation}
        A   = \mathfrak{I}^\top \mathcal{A} \mathfrak{I}, \,
        x         = \mathfrak{I} y, \,
        \tilde{b} = \mathfrak{I}^\top b,
\end{equation}
where $\mathfrak{I}$ is the permutation matrix corresponding to order $I$. The algorithm then proceeds by successively eliminating the sub-blocks associated with the subgraphs in the tree structure, from the leaf level to the top level. 

 An important advantage of nested dissection is both the generation of tree structure and elimination of matrix blocks can be parallellized straightforwardly. However, the effectiveness of nested dissection gradually deteriorates as the elimination progresses to the top of the tree structure. This is because the matrix blocks corresponding to the vertices on the separators become increasingly dense after each elimination. In the end, a dense matrix of order $O(N^{1/2})$ needs to be solved, resulting in an $O(N^{3/2})$ computational complexity. Therefore, to improve the efficiency, an effective solution is required for the dense linear systems on the separators.

\subsection{Sparsification via segments}	
As the key to reduce computational complexity is the efficient factorization of dense blocks on the separators, in this section, we outline a novel LU decomposition technique that recursively sparsifies the intermediate dense blocks. This technique enhances overall efficiency by effectively utilizing geometry information of FEM graph throughout the elimination process, particularly by segments introduced in the following.

\subsubsection{Segments}\label{section_segments}
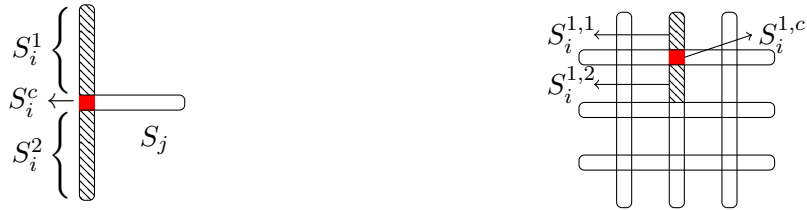
\begin{figure}[tb]
	\centering
    \begin{subfigure} [t] {0.45\textwidth}
        \centering
        \begin{tikzpicture}[
            pattern_slash/.style={pattern=north east lines},
            pattern_Rslash/.style={pattern=north west lines},]

            \filldraw[pattern_Rslash] [rounded corners=2] (-0.1, 1.3) rectangle (0.1, -1.3);
            \fill[red] (-0.1, -0.1) rectangle (0.1, 0.1);
            \draw (-0.1, 0.1) [rounded corners=2] -- (1.3, 0.1) -- (1.3, -0.1) -- (-0.1, -0.1);
            
            \node at (-0.6, 0.7) {$S_i^1 \, \Biggl\{ \Biggr. $};
            \node at (-0.6, -0.7) {$S_i^2 \, \Biggl\{ \Biggr. $};
            \node at (-0.6, 0) {$S_i^c \leftarrow $};
            \node at (0.9, -0.5) {$S_j$};
            
        \end{tikzpicture}
        
        \caption{A separator is split into three segments at level $l$}
        \label{fig-segment (a)}
    \end{subfigure}
    \
    \begin{subfigure} [t] {0.45\textwidth}
        \centering
        \begin{tikzpicture}[
            pattern_slash/.style={pattern=north east lines},
            pattern_Rslash/.style={pattern=north west lines},]

            \fill [pattern_Rslash] (-0.1, 0.1) [rounded corners=2] -- (-0.1, 1.3) -- (0.1, 1.3) -- (0.1, 0.1);
            \draw [rounded corners=2] (-0.1, 1.3) rectangle (0.1, -1.3);
            \fill[red] (-0.1, 0.6) rectangle (0.1, 0.8);
            \draw [rounded corners=2] (1.3, -0.1) rectangle (-1.3, 0.1);

            \draw [rounded corners=2] (1.3, 0.6) rectangle (-1.3, 0.8);
            \draw [rounded corners=2] (1.3, -0.6) rectangle (-1.3, -0.8);
            \draw [rounded corners=2] (0.6, 1.3) rectangle (0.8, -1.3);
            \draw [rounded corners=2] (-0.6, 1.3) rectangle (-0.8, -1.3);

            \node at (-1.4, 1.0) {$S_i^{1,1}$};
            \draw [->] (-0.1, 1.0) -- (-1.1, 1.0);
            \node at (-1.4, 0.34) {$S_i^{1,2}$};
            \draw [->] (-0.1, 0.34) -- (-1.1, 0.34);
            \node at (1.4, 1.0) {$S_i^{1,c}$};
            \draw [->] (0.1, 0.7) -- (1.0, 1.0);
            \node at (0,0) {};
        \end{tikzpicture}
        
        \caption{Segments formed at the lower level $k$ with $k>l$}
        \label{fig-segment (b)}
    \end{subfigure}
	
	\caption{Splitting a separator into segments.
    (a): Separator $S_i$ is split into 3 segments by the separator $S_j$.
    (b): Segment $S_i^1$ is split into smaller segments at the lower level.}
\end{figure}

Motivated by spaND \cite{spaND}, our approach to reduce the complexity of nested dissection focuses on examining the connectivity between different segments. As illustrated in \fig\ref{fig-segment (a)}, during dissection, a separator $S_i$ is divided into three subsets $S_i^1 \cup S_i^2 \cup S_i^c$ by $S_j$, where $S_j$ is a separator at the lower level of the tree structure. Similar to how $S_j$  separates a subgraph of $G$, $S_j \cup S_i^c$ separates $S_i$ into two disjoint sets $S_i^1$ and $S_i^2$. We call $S_i^1$ and $S_i^2$ the \textit{regular} segments of $S_i$, and $S_i^c$  the \textit{junction} segment. In addition, a regular segment like $S_i^1$ can be further subdivided into smaller subsets at the subsequent (or lower) level, as shown in \fig\ref{fig-segment (b)}. These descendant segments, like the subgraphs of $G$, form a hierarchical tree structure within $S_i$.

After eliminating interior vertices at the leaf level, the graph is composed solely of separators, as shown in \fig \ref{fig4}. We refer to this as the \textit{separator graph}, which is associated with the Schur complement of matrix $A$ in \equ(\ref{equ-Ax=b}). Due to the local connectivity property, where edges can not cross any separator, two segments $S_i^k$ and $S_i^l$ from the separator $S_i$, as shown in \fig \ref{fig-finer relation}, are disconnected. Therefore, if we denote by $A_{k,l}$ the submatrix on the separator graph with row indices $k$ associated with  $S_i^k$ and column indices $k$ associated with  $S_i^l$, then the elements in $A_{k,l}$ should be all zero. The same sparsity pattern holds for $A_{l,k}$. This property can be exploited during the elimination process.  

\subsubsection{Sparsification}
\label{section-intro Sparsification}
 In this section, we give a preliminary investigation on the sparsification of the segment $S_i^k\subset S_i$, and its neighbor $V_n$ (i.e., vertex set $V_n = \mathrm{Adj}(S^k_i)$). Let $V_f$ denote the \textit{far}-field of $S^k_i$, where $V_f = \{ v \in V | \ v \ \mathrm{disconnects \ to} \ S^k_i \}$. 
Denote $A_{i,i}$ the submatrix associated with $S^k_i$, $A_{i,n}$ and $A_{n,i}$ the submatrices associated with the interaction between $S^k_i$ and $V_n$,  and $A_{f,n}$ and $A_{n,f}$ the submatrices associated with the interaction between $V_n$ and $V_f$.

Assume that after applying a suitable permutation $E_i$, the Schur complement of $A$, denoted by $A_{s}$, can be written as
\begin{equation}
	E_i^\top A_s E_i = 
	\begin{bmatrix}
		A_{i,i} & A_{i,n} &         \\
		A_{n,i} & A_{n,n} & A_{n,f} \\
		& A_{f,n} & A_{f,f}
	\end{bmatrix}.
	\label{equ-A_{i,f} => 0}
\end{equation}
There exists a transformation matrix $U_{t,i}$, with details provided in Section \ref{sec4}, that factorizes $A_{i,n}$ within a tolerance $\eps$ such that 
\begin{equation}
	U_{t,i}  A_{i,n} = 
	\begin{bmatrix}
		0 \\ 
		A_{s,n}
	\end{bmatrix} + \O(\eps).
\end{equation}
Similarly, the transformation $L_{t,i}$ reduces $A_{n,i}$ to 
$\begin{bmatrix}
	0 & A_{n,s}
\end{bmatrix}$ up to the same tolerance $\eps$.
Define the extended transformation matrices as
$U_{T,i} = \begin{bmatrix}
	U_{t,i} &   \\
	& I
\end{bmatrix}$
and
$L_{T,i} = \begin{bmatrix}
	L_{t,i} &   \\
	& I
\end{bmatrix}$ with $I$ the identity matrix such that $U_{T,i} E_i^\top A E_i L_{T,i}$ is well-defined. Then we have
\begin{equation}
	U_{T,i} E_i^\top A E_i L_{T,i} = 
	\begin{bmatrix}
		A_{r,r} & A_{r,s} & 0 &  \\
		A_{s,r} & A_{s,s} & A_{s,n} & \\
		0               & A_{n,s} & A_{n,n}         & A_{n,f} \\
		                &                 & A_{f,n}         & A_{f,f}
	\end{bmatrix} + \O(\eps).
	\label{equ-UtALt}
\end{equation}

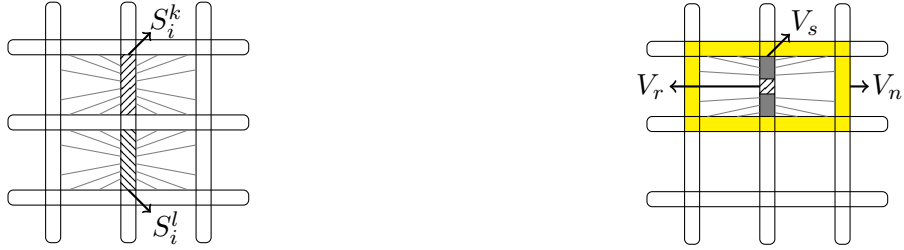
\begin{figure}[tb]
	\centering
	\begin{subfigure} [t] {0.48\textwidth}
		\centering
        \begin{tikzpicture}[scale=1, 
            pattern_slash/.style={pattern=north east lines},
            pattern_Rslash/.style={pattern=north west lines},]
            
            \fill[pattern_Rslash] (0.1,1.1) rectangle (-0.1,1.9);
            \fill[pattern_slash] (0.1,2.1) rectangle (-0.1,2.9);
            \draw [rounded corners=2] (-0.1,0.4) rectangle (0.1,3.6);
            
            \draw [rounded corners=2] (-1.6,0.9) rectangle (1.6,1.1);
            \draw [rounded corners=2] (-1.6,1.9) rectangle (1.6,2.1);
            \draw [rounded corners=2] (-1.6,2.9) rectangle (1.6,3.1);
            
            \draw [rounded corners=2] (-0.9,0.4) rectangle (-1.1,3.6);
            \draw [rounded corners=2] (0.9,0.4) rectangle (1.1,3.6);
            
            \foreach \j in {0, -1}
            {
                \foreach \i in {1, -1}
                {
                    \draw [gray] (\i*0.1,2.25+\j) -- (\i*0.4,2.1+\j);
                    \draw [gray] (\i*0.1,2.35+\j) -- (\i*0.8,2.1+\j);
                    \draw [gray] (\i*0.1,2.45+\j) -- (\i*0.9,2.3+\j);
                    \draw [gray] (\i*0.1,2.55+\j) -- (\i*0.9,2.7+\j);
                    \draw [gray] (\i*0.1,2.65+\j) -- (\i*0.8,2.9+\j);
                    \draw [gray] (\i*0.1,2.75+\j) -- (\i*0.4,2.9+\j);
                }
            }

            \node at (0.5,3.35) {$S_i^k$};
            \draw [thick] [->] (0,2.9) -- (0.3, 3.2);
            \node at (0.5,0.6) {$S_i^l$};
            \draw [thick] [->] (0,1.1) -- (0.3, 0.8);
        \end{tikzpicture}
		
		\caption{Edge relations for two segments from one separator}
		\label{fig-finer relation}
	\end{subfigure}
	\quad
	\begin{subfigure} [t] {0.48\textwidth}
		\centering
		\begin{tikzpicture}[scale=1, 
            pattern_slash/.style={pattern=north east lines},
            pattern_Rslash/.style={pattern=north west lines},]
            
            \fill[pattern_slash] (0.1,2.4) rectangle (-0.1,2.6);
            \fill[gray] (0.1,2.6) rectangle (-0.1,2.9);
            \draw (-0.1,2.6) -- (0.1,2.6);
            \fill[gray] (0.1,2.1) rectangle (-0.1,2.4);
            \draw (-0.1,2.4) -- (0.1,2.4);

            \fill[yellow] (-1.1,1.9) rectangle (1.1,2.1);
            \fill[yellow] (-1.1,2.9) rectangle (1.1,3.1);
            \fill[yellow] (-1.1,2.1) rectangle (-0.9,2.9);
            \fill[yellow] (1.1,2.1) rectangle (0.9,2.9);
            
            \draw [rounded corners=2] (-0.1,0.4) rectangle (0.1,3.6);
            
            \draw [rounded corners=2] (-1.6,0.9) rectangle (1.6,1.1);
            \draw [rounded corners=2] (-1.6,1.9) rectangle (1.6,2.1);
            \draw [rounded corners=2] (-1.6,2.9) rectangle (1.6,3.1);
            
            \draw [rounded corners=2] (-0.9,0.4) rectangle (-1.1,3.6);
            \draw [rounded corners=2] (0.9,0.4) rectangle (1.1,3.6);
            
            \foreach \i in {1, -1}
            {
                \draw [gray] (\i*0.1,2.15) -- (\i*0.4,2.1);
                \draw [gray] (\i*0.1,2.25) -- (\i*0.8,2.1);
                \draw [gray] (\i*0.1,2.35) -- (\i*0.9,2.3);
                \draw [gray] (\i*0.1,2.65) -- (\i*0.9,2.7);
                \draw [gray] (\i*0.1,2.75) -- (\i*0.8,2.9);
                \draw [gray] (\i*0.1,2.85) -- (\i*0.4,2.9);
            }
            
            \node at (0,0.4) {};

            \node at (0.5,3.35) {$V_s$};
            \draw [thick] [->] (0,2.9) -- (0.3, 3.2);
            \node at (-1.55,2.5) {$V_r$};
            \draw [thick] [->] (-0.1, 2.5) -- (-1.3, 2.5);
            \node at (1.6,2.5) {$V_n$};
            \draw [thick] [->] (1.1, 2.5) -- (1.35, 2.5);
            
        \end{tikzpicture}
		
		\caption{Segments after sparsification}
		\label{fig-sparse example}
	\end{subfigure}
	
	\caption{Sparsification of segments. (a): Two segments $S_i^l$ and $S_i^k$ belong to the same separator $S_i$ but have disjoint edge relations.
    (b): After the sparsification, the skeleton $V_s$ of $S_i^k$ still connects to $V_n$, the neighbor of $S_i^k$, but its remainder $V_r$ does not connect to $V_n$ anymore.}\label{fig4}
\end{figure}

From \equ(\ref{equ-UtALt}), it becomes clear that vertex set $V_s$, corresponding to the row indices of $A_{s,n}$, remains connected to the neighboring vertices of $S_i^k$, that is, the vertex set $V_n$. However, $V_r = S^k_i \backslash V_s$ is no longer connected to $V_n$, indicating that $V_s$ can represent the entire segment $S^k_i$ in interaction with other segments. Consequently, the corresponding matrix elements of $V_r$ is decoupled from $V_n$ (i.e., $A_{r,n} = 0$, $A_{n,r} = 0$), 
as illustrated in \fig \ref{fig-sparse example}. We refer to $V_s$ as the \textit{skeleton} of $S^k_i$ and $V_r$ the \textit{remainder} of $S^k_i$.
This represents a single step in the sparsification of matrix $A$, where the transformations $U_{T,i}$ and $L_{T,i}$ eliminate the non-zeros elements associated with remainder $V_r$ without introducing extra \textit{fill-in}s in the rest of matrix $A$. In particular, $A_{f,f}$, $A_{n,f}$ and $A_{f,n}$ remain unchanged, so the sparsification is completely local. Once all segments at level $l$ are sparsified, the corresponding  skeletons can be merged to form new segments at the next higher level. 

Especially, we combine all the transformations at level $l$ to form two block-diagonal matrices,  $\mathcal{L}_l$ and $\mathcal{U}_l$, with each block corresponding to the sparsification of a single segment. These two matrices have the forms:
\begin{equation}
\begin{aligned}
\mathcal{L}_l &= \prod_{i=I_l}^1 L_{T,i,l}^{-1} E_{i,l}^\top, \quad
\mathcal{U}_l &= \prod_{i=1}^{I_l} E_{i,l} U_{T,i,l}^{-1},
\end{aligned}
\end{equation}
where $\{1, 2, \dots I_l\}$ represents all the non-junction segments at level $l$, $L_{T,i,l}$ and $U_{T,i,l}$ are the transformation matrices for the $i$-th segment, and $E_{i,l}$ is the permutation matrix for the $i$-th segment at level $l$. After merging the skeletons, we can recursively apply these transformations to each segment in the higher level. In the end, it results in a decomposition of matrix $A$:
\begin{equation}
A \approx \mathcal{U}_L \mathcal{U}_{L-1}\cdots \mathcal{U}_1 \mathcal{L}_1 \cdots \mathcal{L}_{L-1} \mathcal{L}_L.
\label{equ-spLU_lite}
\end{equation}
We call the factorization \eqref{equ-spLU_lite} as the recursive sparse LU decomposition, or spaLU, of $A$. Further details will be provided in Section \ref{sec4}.

\section{Nested structure for finite element graph}
In this section, we give the details of how to construct the necessary data structures for separators and segments from a finite element graph for use in the recursive LU decomposition. 

\subsection{Structure of separators}
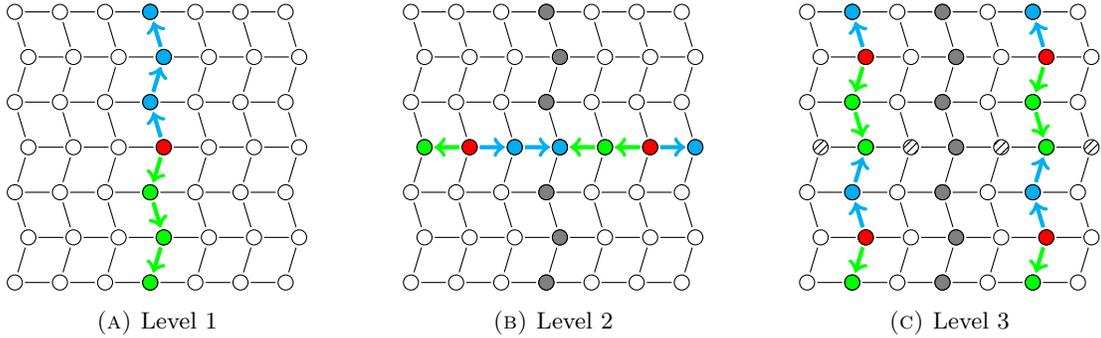
\begin{figure}[tb]
	\centering
	\begin{subfigure} [t] {0.3\textwidth}
        \centering
		\begin{tikzpicture} [scale = 0.6]
			\tikzmath{ \r = 0.1 / 0.6; \s = 0.3; };
			\foreach \j in {-3, -1, 1, 3}
			{
				\foreach \i in {-3, -2, -1, 0, 1, 2, 3}
				{
					\node (a\j\i) at (\i, \j) {};
					\filldraw [fill=white] (a\j\i) circle [radius=\r];
				}
			}
            \foreach \j in {-2, 0, 2}
			{
				\foreach \i in {-3, -2, -1, 0, 1, 2, 3}
				{
					\node (a\j\i) at (\i + \s, \j) {};
					\filldraw [fill=white] (a\j\i) circle [radius=\r];
				}
			}
            \foreach \j in {-3, -2, -1, 0, 1, 2, 3}
            {
                \foreach \i/\k in {-3/-2, -2/-1, -1/0, 0/1, 1/2, 2/3}
                {
                    \draw (a\i\j) -- (a\k\j);
                }
            }
            \foreach \j in {-3, -2, -1, 0, 1, 2, 3}
            {
                \foreach \i/\k in {-3/-2, -2/-1, -1/0, 0/1, 1/2, 2/3}
                {
                    \draw (a\j\i) -- (a\j\k);
                }
            }
			
			\fill[fill=red,draw=black] (a00) circle [radius=\r];
			\foreach \i/\j in {1/0, 2/1, 3/2}
			{
				\fill[fill=cyan,draw=black] (a\i0) circle [radius=\r];
				\draw[->][ultra thick][cyan] (a\j0) -- (a\i0);
			}
			\foreach \i/\j in {-1/0, -2/-1, -3/-2}
			{
				\fill[fill=green,draw=black] (a\i0) circle [radius=\r];
				\draw[->][ultra thick][green] (a\j0) -- (a\i0);
			}
		\end{tikzpicture}
    \caption{Level 1}
	\end{subfigure}
	\
	\begin{subfigure} [t] {0.3\textwidth}
        \centering
		\begin{tikzpicture} [scale = 0.6]
			\tikzmath{ \r = 0.1 / 0.6; \s = 0.3; };
			\foreach \j in {-3, -1, 1, 3}
			{
				\foreach \i in {-3, -2, -1, 0, 1, 2, 3}
				{
					\node (a\j\i) at (\i, \j) {};
					\filldraw [fill=white] (a\j\i) circle [radius=\r];
				}
			}
            \foreach \j in {-2, 0, 2}
			{
				\foreach \i in {-3, -2, -1, 0, 1, 2, 3}
				{
					\node (a\j\i) at (\i + \s, \j) {};
					\filldraw [fill=white] (a\j\i) circle [radius=\r];
				}
			}
            \foreach \j in {-3, -2, -1, 0, 1, 2, 3}
            {
                \foreach \i/\k in {-3/-2, -2/-1, -1/0, 0/1, 1/2, 2/3}
                {
                    \draw (a\i\j) -- (a\k\j);
                }
            }
            \foreach \j in {-3, -2, -1, 0, 1, 2, 3}
            {
                \foreach \i/\k in {-3/-2, -2/-1, -1/0, 0/1, 1/2, 2/3}
                {
                    \draw (a\j\i) -- (a\j\k);
                }
            }
			
			\foreach \i in {-3, -2, -1, 0, 1, 2, 3}
			{
				\filldraw [fill=gray] (a\i0) circle [radius=\r];
			}
			
			\fill[fill=red,draw=black] (a02) circle [radius=\r];
			\fill[fill=red,draw=black] (a0-2) circle [radius=\r];
			\foreach \i/\j in {3/2}
			{
				\fill[fill=cyan,draw=black] (a0\i) circle [radius=\r];
				\draw[->][ultra thick][cyan] (a0\j) -- (a0\i);
			}
			\foreach \i/\j in {1/2, 0/1}
			{
				\fill[fill=green,draw=black] (a0\i) circle [radius=\r];
				\draw[->][ultra thick][green] (a0\j) -- (a0\i);
			}
			\foreach \i/\j in {-3/-2}
			{
				\fill[fill=green,draw=black] (a0\i) circle [radius=\r];
				\draw[->][ultra thick][green] (a0\j) -- (a0\i);
			}
			\foreach \i/\j in {-1/-2, 0/-1}
			{
				\fill[fill=cyan,draw=black] (a0\i) circle [radius=\r];
				\draw[->][ultra thick][cyan] (a0\j) -- (a0\i);
			}
		\end{tikzpicture}
    \caption{Level 2}
	\end{subfigure}
	\
	\begin{subfigure} [t] {0.3\textwidth}
        \centering
		\begin{tikzpicture} [scale = 0.6,
			pattern_slash/.style={pattern=north east lines},]
			\tikzmath{ \r = 0.1 / 0.6; \s = 0.3; };
			\foreach \j in {-3, -1, 1, 3}
			{
				\foreach \i in {-3, -2, -1, 0, 1, 2, 3}
				{
					\node (a\j\i) at (\i, \j) {};
					\filldraw [fill=white] (a\j\i) circle [radius=\r];
				}
			}
            \foreach \j in {-2, 0, 2}
			{
				\foreach \i in {-3, -2, -1, 0, 1, 2, 3}
				{
					\node (a\j\i) at (\i + \s, \j) {};
					\filldraw [fill=white] (a\j\i) circle [radius=\r];
				}
			}
            \foreach \j in {-3, -2, -1, 0, 1, 2, 3}
            {
                \foreach \i/\k in {-3/-2, -2/-1, -1/0, 0/1, 1/2, 2/3}
                {
                    \draw (a\i\j) -- (a\k\j);
                }
            }
            \foreach \j in {-3, -2, -1, 0, 1, 2, 3}
            {
                \foreach \i/\k in {-3/-2, -2/-1, -1/0, 0/1, 1/2, 2/3}
                {
                    \draw (a\j\i) -- (a\j\k);
                }
            }
			
			\foreach \i in {-3, -2, -1, 0, 1, 2, 3}
			{
				\filldraw [fill=gray] (a\i0) circle [radius=\r];
			}
			\foreach \i in {1, 2, 3}
			{
				\fill [pattern_slash] (a0\i) circle [radius=\r];
				\fill [pattern_slash] (a0-\i) circle [radius=\r];
			}
			
			\fill[fill=red,draw=black] (a22) circle [radius=\r];
			\fill[fill=red,draw=black] (a-22) circle [radius=\r];
			\fill[fill=red,draw=black] (a2-2) circle [radius=\r];
			\fill[fill=red,draw=black] (a-2-2) circle [radius=\r];
			
			\foreach \i/\j in {3/2, -1/-2, 0/-1}
			{
				\foreach \k in {2, -2}
				{
					\fill[fill=cyan,draw=black] (a\i\k) circle [radius=\r];
					\draw[->][ultra thick][cyan] (a\j\k) -- (a\i\k);
				}
			}
			\foreach \i/\j in {1/2, 0/1, -3/-2}
			{
				\foreach \k in {2, -2}
				{
					\fill[fill=green,draw=black] (a\i\k) circle [radius=\r];
					\draw[->][ultra thick][green] (a\j\k) -- (a\i\k);
				}
			}
		\end{tikzpicture}
    \caption{Level 3}
	\end{subfigure}
	
	\caption {A sketch of generating separators. Given the separator direction $\vec{d}$ in each subgraph,  the separator is generated by expanding from the center vertex (in red) along the  direction $\vec{d}$ (in blue) and $-\vec{d}$ (in green). }
	\label{fig-separator 2D}
\end{figure}

For a two dimensional finite element graph, since vertices may not be uniformly distributed, as shown in Figure \ref{fig-separator 2D}, the key strategy for separating the graph is to identify a polygonal line composed of edges that expands from the graph's center vertex $c$. Given a graph $G = (V,E)$, we begin by estimating the width and height of the bounding box that contains $V$. The center of $V$, denoted as $v_c$, is determined by calculating the median of the coordinates. The vertex $c\in V$  closest to $v_c$ is then selected as the center vertex. The unit vector $\vec{d}$ along the positive x-axis is chosen as the separating direction if the width of $V$ is less than the height; otherwise, the positive y-axis is selected as the separating direction.

Define the separator generating function
\begin{equation}\label{degreebias}
	D(u,v,c,\vec{d}) = \frac{\left< u-v,\vec{d} \right>}{\norm(u-v) \norm(\vec{d})} + 
	\theta \frac{\left< u-c,\vec{d} \right>}{\norm(u-c) \norm(\vec{d})}, 
\end{equation}
where $\left< \cdot, \cdot \right>$ denotes the inner product in $\mathbb{R}^2$. This function $D(u,v,c,\vec{d})$ provides a degree bias relative to the direction $\vec{d}$. The first term in $D$ measures the alignment of $u-v$ with $\vec{d}$, while the second term, scaled by a parameter $\theta$ ($0 < \theta < 1$), measures the deviation of the alignment of $u-c$ from $\vec{d}$. 

Here we generate the separator by choosing $u=Next(v,\vec{d})$ that maximizes $D(u,v,c,\vec{d})$. Specifically, we take $v_0 = c$, $\theta=0.1$, and iteratively choose $v_{i+1} = Next(v_i,\vec{d})$, $i = 1,\dots,p$, until $v_{p+1} = Next(v_p,\vec{d})$ does not exist in $V$ or $D(v_{p+1},v_{p},c,\vec{d}) \le 0$. Similarly, vertices $v_{j-1} = Next(v_j, -\vec{d})$, $j = -1,\dots,-q$, are obtained until $v_{-q-1}$ does not exist in $V$ or $D(v_{-q-1},v_{-q},c,-\vec{d}) \le 0$. The separator $S$ is then composed of the sequence $v_{-q},\dots,v_{-1},v_0,v_1,\dots,v_p$, which divides the graph $G$ into two subgraphs $G_1$ and $G_2$ by construction. By recursively applying this procedure, one can identify separators for the subgraphs $G_1$ and $G_2$ and their successors, as illustrated in \fig\ref{fig-separator 2D}. 

\begin{algorithm}[htb]
	\caption{Nested dissection in 2D with segment information.}
	\label{alg-nested 2D}
	\begin{algorithmic}[1]
		\Require A set of $n$ vertices $V$ with their coordinates, and the associated sparse matrix $\mathcal{A}$.
		\Ensure Nested dissection order $I$, and segment information $\mathcal{N}$.
		\State Extract edge relations $E$ from the structure of matrix $\mathcal{A}$ and let the number of layers $L =  \O(\log(N))$. (Here $\log$ is logarithm with base 2.)
		\State Let $B$ be the set of boundary segments of $V$ and $B = \emptyset$. \label{alg-nested 2D - step 2}
        \State Initialize $I \gets \emptyset$, $l \gets 1$ and $\mathcal{N} \gets \{\mathcal{N}_1, \dots, \mathcal{N}_L\} = \{ \emptyset, \dots, \emptyset \}$
		
		\Function{Nested2D}{$V, B, l$}
		\State Find the center vertex $c$ of $V$ and the separating direction $\vec{d}$.
        \State Set the separator $S \gets \{ c \}$.
		\State Append vertices from $V$ to $S$ by expanding from $c$ along $\vec{d}$  and  $-\vec{d}$ using equation \eqref{degreebias}.
		\State Update $I \gets I \cup S$.
		\State Partition $V$ into two parts $V_1$ and $V_2$ by the separator $S$.
        \For {each boundary segment $B^k \in B$}
        \State Determine the junction $S' \gets$ intersection of $S$ and $B^k$.
            \If {$S' \ne \emptyset$}
                \State  Split $B^k$ into $B^k_1$ and $B^k_2$, and replace $B^k$ by $B^k_1$ and $B^k_2$ in $B$ and $\mathcal{N}_l$.
            \EndIf
        \EndFor
        \State Set $B_1 \gets \{ B^{k} \in B \ |\  B^{k} \text{ on one side of } S \} \cup \{S\}$
        \State Set $B_2 \gets \{ B^{k} \in B \ |\  B^{k} \text{ on the other side of } S \} \cup \{S\}$.
        \State Update $\mathcal{N}_l \gets \{ \mathcal{N}_l, S \}$.
		\For{$i=1,2$}
    		\If{$l < L$}
    		\State \Call{Nested2D}{$V_i, B_i, l+1$}.
    		\Else
    		\State $I \gets I \cup V_i$.
    		\EndIf
		\EndFor
		
		\State \Return $I, \mathcal{N}$.
		\EndFunction
	\end{algorithmic}
\end{algorithm}


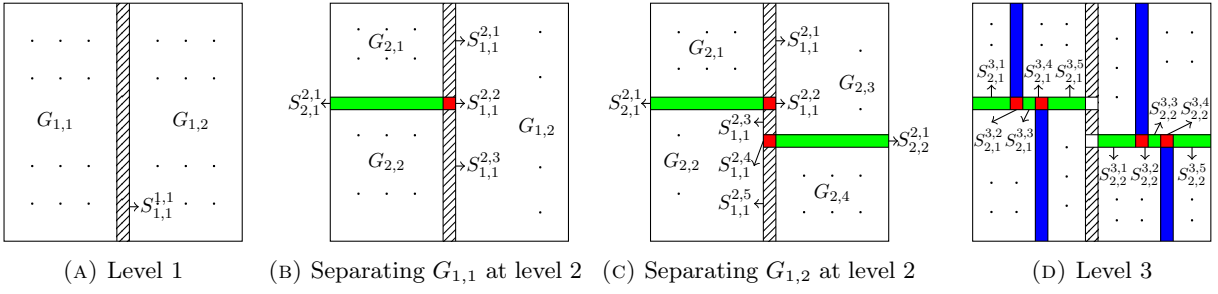
\begin{figure}[tb]
    \centering
    \begin{subfigure} [t] {0.2\textwidth}
		\centering
		\begin{tikzpicture}[scale=0.8333, every node/.append style={scale=0.7},
			pattern_slash/.style={pattern=north east lines},
			pattern_Rslash/.style={pattern=north west lines},]
			\tikzmath{ \b = 1.8; };
			
			\fill [pattern_slash] (-0.1, -\b - 0.1) rectangle (0.1, \b + 0.1);
			
			\draw (-\b - 0.1, -\b - 0.1) rectangle (\b + 0.1, \b + 0.1);
			\draw (-0.1, -\b - 0.1) -- (-0.1, \b + 0.1);
			\draw (0.1, -\b - 0.1) -- (0.1, \b + 0.1);
			
			\foreach \i in {1, 2, 3}
			{
				\foreach \j in {1, 2}
				{
					\filldraw (\b / 4 * \i + 0.1, \b / 3 * \j + 0.1) circle [radius = 0.01];
					\filldraw (\b / 4 * \i + 0.1, -\b / 3 * \j - 0.1) circle [radius = 0.01];
					\filldraw (-\b / 4 * \i - 0.1, \b / 3 * \j + 0.1) circle [radius = 0.01];
					\filldraw (-\b / 4 * \i - 0.1, -\b / 3 * \j - 0.1) circle [radius = 0.01];
				}
			}
			\node at (-0.6 * \b, 0) {$G_{1,1}$};
			\node at (0.6 * \b, 0) {$G_{1,2}$};
			\draw [->] (0.1, -\b * 0.75) -- (0.25, -\b * 0.75);
			\node at (0.55, -\b * 0.75) {$S_{1,1}^{1,1}$};
			
		\end{tikzpicture}
		
		\caption{Level 1}
	\end{subfigure}
	\
	\begin{subfigure} [t] {0.25\textwidth}
		\centering
		\begin{tikzpicture}[scale=0.8333, every node/.append style={scale=0.7},
			pattern_slash/.style={pattern=north east lines},
			pattern_Rslash/.style={pattern=north west lines},]
			\tikzmath{
				\b = 1.8;
				\len1 = (\b + 0.3) * 2/9;
				\len2 = (\b + 0.3) * 4/9;
				\len3 = \b * 2/5;
			};
			
			\fill [pattern_slash] (-0.1, -\b-0.1) rectangle (0.1, \b+0.1);
			\fill [red] (-0.1, \b / 6 - 0.1) rectangle (0.1, \b / 6 + 0.1);
			\fill [green] (-\b - 0.1, \b / 6 - 0.1) rectangle (-0.1, \b / 6 + 0.1);
			
			\draw (-\b - 0.1, -\b - 0.1) rectangle (\b + 0.1, \b + 0.1);
			\draw (-0.1, -\b - 0.1) -- (-0.1, \b + 0.1);
			\draw (0.1, -\b - 0.1) -- (0.1, \b + 0.1);
			\draw (-\b - 0.1, \b / 6 - 0.1) -- (0.1, \b / 6 - 0.1);
			\draw (-\b - 0.1, \b / 6 + 0.1) -- (0.1, \b / 6 + 0.1);
			
			\foreach \i in {1, 2, 3}
			{
				\foreach \j in {1, 2}
				{
					\filldraw (-\b / 4 * \i - 0.1, \len1 * \j + \b/3 - 0.05) circle [radius = 0.01];
					\filldraw (-\b / 4 * \i - 0.1, -\len2 * \j + \b/3 + 0.125) circle [radius = 0.01];
				}
			}
			\foreach \j in {-2, -1, 1, 2}
			{
				\filldraw (\b * 3/4 + 0.1, \len3 * \j) circle [radius = 0.01];
			}
			
			\node at (-0.5 * \b - 0.1, \b * 0.75 - 0.1) {$G_{2,1}$};
			\node at (-0.5 * \b - 0.1, -\b * 0.25 - 0.2) {$G_{2,2}$};
			\node at (0.5 * \b + 0.5, -0.1) {$G_{1,2}$};
			\draw [->](-\b - 0.1, \b / 6) -- (-\b - 0.25, \b / 6);
			\node at (-\b - 0.5, \b / 6) {$S_{2,1}^{2,1}$};
			
			\draw [->] (0.1, \b /3*2 + 0.1) -- (0.25, \b /3*2 + 0.1);
			\node at (0.55, \b /3*2 + 0.1) {$S_{1,1}^{2,1}$};
			\draw [->] (0.1, \b / 6) -- (0.25, \b / 6);
			\node at (0.55, \b / 6) {$S_{1,1}^{2,2}$};
			\draw [->] (0.1, -\b /3 - 0.1) -- (0.25, -\b /3 - 0.1);
			\node at (0.55, -\b /3 - 0.1) {$S_{1,1}^{2,3}$};

		\end{tikzpicture}
		
		\caption{Separating $G_{1,1}$ at level 2}
	\end{subfigure}
	\
	\begin{subfigure} [t] {0.25\textwidth}
		\centering
		\begin{tikzpicture}[scale=0.8333, every node/.append style={scale=0.7},
			pattern_slash/.style={pattern=north east lines},
			pattern_Rslash/.style={pattern=north west lines},]
			\tikzmath{
				\b = 1.8;
				\len1 = (\b + 0.3) * 5/18;
				\len2 = (\b + 0.3) * 4/9;
			};
			
			\fill [pattern_slash] (-0.1, -\b-0.1) rectangle (0.1, \b+0.1);
			\fill [red] (-0.1, \b / 6 - 0.1) rectangle (0.1, \b / 6 + 0.1);
			\fill [green] (-\b - 0.1, \b / 6 - 0.1) rectangle (-0.1, \b / 6 + 0.1);
			\fill [red] (-0.1, -\b / 6 - 0.1) rectangle (0.1, -\b / 6 + 0.1);
			\fill [green] (\b + 0.1, -\b / 6 - 0.1) rectangle (0.1, -\b / 6 + 0.1);
			
			\draw (-\b - 0.1, -\b - 0.1) rectangle (\b + 0.1, \b + 0.1);
			\draw (-0.1, -\b - 0.1) -- (-0.1, \b + 0.1);
			\draw (0.1, -\b - 0.1) -- (0.1, \b + 0.1);
			\draw (-\b - 0.1, \b / 6 - 0.1) -- (0.1, \b / 6 - 0.1);
			\draw (-\b - 0.1, \b / 6 + 0.1) -- (0.1, \b / 6 + 0.1);
			\draw (\b + 0.1, -\b / 6 - 0.1) -- (-0.1, -\b / 6 - 0.1);
			\draw (\b + 0.1, -\b / 6 + 0.1) -- (-0.1, -\b / 6 + 0.1);

			\foreach \i in {1, 2, 3}
			{
				\foreach \j in {1, 2}
				{
					\filldraw (-\b / 4 * \i - 0.1, \len1 * \j + \b/6 - 0.03) circle [radius = 0.01];
					\filldraw (\b / 4 * \i + 0.1, -\len1 * \j - \b/6 + 0.03) circle [radius = 0.01];
				}
			}
			\foreach \j in {1, 2}
			{
				\filldraw (\b * 3/4 + 0.1, \len2 * \j - \b/3 - 0.125) circle [radius = 0.01];
				\filldraw (-\b * 3/4 - 0.1, -\len2 * \j + \b/3 + 0.125) circle [radius = 0.01];
			}
			
			\node at (-0.5 * \b - 0.1, \b * 7/12 + 0.1) {$G_{2,1}$};
			\node at (-0.5 * \b - 0.5, -\b * 0.25 - 0.2) {$G_{2,2}$};
			\node at (0.5 * \b + 0.5, \b / 3) {$G_{2,3}$};
			\node at (0.5 * \b + 0.1, -\b * 7/12 - 0.1) {$G_{2,4}$};
			\draw [->](-\b - 0.1, \b /6) -- (-\b - 0.25, \b /6);
			\node at (-\b - 0.5, \b /6) {$S_{2,1}^{2,1}$};
			\draw [->](\b + 0.1, -\b /6) -- (\b + 0.25, -\b /6);
			\node at (\b + 0.5, -\b /6) {$S_{2,2}^{2,1}$};
			
			\draw [->] (0.1, \b /3*2 + 0.1) -- (0.25, \b /3*2 + 0.1);
			\node at (0.55, \b /3*2 + 0.1) {$S_{1,1}^{2,1}$};
			\draw [->] (0.1, \b / 6) -- (0.25, \b / 6);
			\node at (0.55, \b / 6) {$S_{1,1}^{2,2}$};
			\draw [->] (-0.1, 0) -- (-0.25, 0);
			\node at (-0.55, 0 - 0.1) {$S_{1,1}^{2,3}$};
			\draw [->] (-0.1, -\b / 6) -- (-0.25, -\b / 3 - 0.1);
			\node at (-0.55, -\b / 3 - 0.1) {$S_{1,1}^{2,4}$};
			\draw [->] (-0.1, -\b / 3*2 - 0.1) -- (-0.25, -\b / 3*2 - 0.1);
			\node at (-0.55, -\b / 3*2 - 0.1) {$S_{1,1}^{2,5}$};
			
		\end{tikzpicture}
		
		\caption{Separating $G_{1,2}$ at level 2}
	\end{subfigure}
	\
	\begin{subfigure} [t] {0.25\textwidth}
		\centering
		\begin{tikzpicture}[scale=0.8333, every node/.append style={scale=0.6},
			pattern_slash/.style={pattern=north east lines},
			pattern_Rslash/.style={pattern=north west lines},]
			\tikzmath{
				\b = 1.8;
				\a1 = \b * 3/4 + 0.25;
				\a2 = \b / 2 + 0.3;
				\a3 = \b / 2 - 0.1;
				\a4 = \b / 4 - 0.1;
				\len1 = (\b - 0.6) * 5/6 / 3;
				\len2 = (\b - 0.6) * 7/6 / 3;
			};
			
			\fill [pattern_slash] (-0.1, -\b-0.1) rectangle (0.1, \b+0.1);
			\fill [white] (-0.1, \b / 6 - 0.1) rectangle (0.1, \b / 6 + 0.1);
			\fill [green] (-\b - 0.1, \b / 6 - 0.1) rectangle (-0.1, \b / 6 + 0.1);
			\fill [white] (-0.1, -\b / 6 - 0.1) rectangle (0.1, -\b / 6 + 0.1);
			\fill [green] (\b + 0.1, -\b / 6 - 0.1) rectangle (0.1, -\b / 6 + 0.1);
			
			\fill [blue] (-\b / 2 - 0.2, \b / 6 + 0.1) rectangle (-\b / 2 - 0.4, \b + 0.1);
			\fill [blue] (\b / 2 - 0.2, -\b / 6 + 0.1) rectangle (\b / 2, \b + 0.1);
			\fill [blue] (-\b / 2 + 0.2, \b / 6 - 0.1) rectangle (-\b / 2, -\b - 0.1);
			\fill [blue] (\b / 2 + 0.2, -\b / 6 - 0.1) rectangle (\b / 2 + 0.4, -\b - 0.1);
			\fill [red] (-\b / 2 - 0.4, \b / 6 + 0.1) rectangle (-\b / 2 - 0.2, \b / 6 - 0.1);
			\fill [red] (-\b / 2, \b / 6 + 0.1) rectangle (-\b / 2 + 0.2, \b / 6 - 0.1);
			\fill [red] (\b / 2 - 0.2, -\b / 6 + 0.1) rectangle (\b / 2, -\b / 6 - 0.1);
			\fill [red] (\b / 2 + 0.2, -\b / 6 + 0.1) rectangle (\b / 2 + 0.4, -\b / 6 - 0.1);
			
			\draw (-\b - 0.1, -\b - 0.1) rectangle (\b + 0.1, \b + 0.1);
			\draw (-0.1, -\b - 0.1) -- (-0.1, \b + 0.1);
			\draw (0.1, -\b - 0.1) -- (0.1, \b + 0.1);
			\draw (-\b - 0.1, \b / 6 - 0.1) -- (0.1, \b / 6 - 0.1);
			\draw (-\b - 0.1, \b / 6 + 0.1) -- (0.1, \b / 6 + 0.1);
			\draw (\b + 0.1, -\b / 6 - 0.1) -- (-0.1, -\b / 6 - 0.1);
			\draw (\b + 0.1, -\b / 6 + 0.1) -- (-0.1, -\b / 6 + 0.1);
			
			\draw (-\b / 2 - 0.4, \b / 6 - 0.1) -- (-\b / 2 - 0.4, \b + 0.1);
			\draw (-\b / 2 - 0.2, \b / 6 - 0.1) -- (-\b / 2 - 0.2, \b + 0.1);
			\draw (-\b / 2, \b / 6 + 0.1) -- (-\b / 2, -\b - 0.1);
			\draw (-\b / 2 + 0.2, \b / 6 + 0.1) -- (-\b / 2 + 0.2, -\b - 0.1);
			
			\draw (\b / 2 - 0.2, -\b / 6 - 0.1) -- (\b / 2 - 0.2, \b + 0.1);
			\draw (\b / 2, -\b / 6 - 0.1) -- (\b / 2, \b + 0.1);
			\draw (\b / 2 + 0.2, -\b / 6 + 0.1) -- (\b / 2 + 0.2, -\b - 0.1);
			\draw (\b / 2 + 0.4, -\b / 6 + 0.1) -- (\b / 2 + 0.4, -\b - 0.1);
			
			\foreach \i in {1, 2}
			{
				\filldraw (-\a1, \b / 6 + 0.6 + \len1 * \i) circle [radius = 0.01];
				
				\filldraw (-\a3, \b / 6 + 0.6 + \len1 * \i) circle [radius = 0.01];
				\filldraw (-\a4 - 0.05, \b / 6 + 0.6 + \len1 * \i) circle [radius = 0.01];
				
				\filldraw (-\a1, -\b / 6 - 0.2 - \len2 * \i) circle [radius = 0.01];
				\filldraw (-\a2, -\b / 6 - 0.2 - \len2 * \i) circle [radius = 0.01];
				
				\filldraw (\a1, \b / 6 + 0.2 + \len2 * \i) circle [radius = 0.01];
				\filldraw (\a2, \b / 6 + 0.2 + \len2 * \i) circle [radius = 0.01];
				
				\filldraw (\a3, -\b / 6 - 0.6 - \len1 * \i) circle [radius = 0.01];
				\filldraw (\a4 + 0.05, -\b / 6 - 0.6 - \len1 * \i) circle [radius = 0.01];
				
				\filldraw (\a1, -\b / 6 - 0.6 - \len1 * \i) circle [radius = 0.01];
			}
			\foreach \i in {0, 1, 2}
			{
				\filldraw (-\a4 - 0.05, -\b / 6 - 0.1 - \len2 * \i) circle [radius = 0.01];
				\filldraw (\a4 + 0.05, \b / 6 + 0.1 + \len2 * \i) circle [radius = 0.01];
			}
			
			\node at (-\a1, \b / 6 + 0.5) {$S_{2,1}^{3,1}$};
			\draw [->] (-\a1, \b / 6 + 0.1) -- (-\a1, \b / 6 + 0.3);
			\node at (-\a1 - 0.05, \b / 6 - 0.6) {$S_{2,1}^{3,2}$};
			\draw [->] (-\a2, \b / 6 - 0.1) -- (-\a1, \b / 6 - 0.3);
			\node at (-\a2 + 0.05, \b / 6 - 0.6) {$S_{2,1}^{3,3}$};
			\draw [->] (-\a2 + 0.2, \b / 6 - 0.1) -- (-\a2 + 0.1, \b / 6 - 0.3);
			\node at (-\a3 - 0.05, \b / 6 + 0.5) {$S_{2,1}^{3,4}$};
			\draw [->] (-\a3 - 0.05, \b / 6 + 0.1) -- (-\a3 - 0.05, \b / 6 + 0.3);
			\node at (-\a4, \b / 6 + 0.5) {$S_{2,1}^{3,5}$};
			\draw [->] (-\a4, \b / 6 + 0.1) -- (-\a4, \b / 6 + 0.3);
			
			\node at (\a1, -\b / 6 - 0.55) {$S_{2,2}^{3,5}$};
			\draw [->] (\a1, -\b / 6 - 0.1) -- (\a1, -\b / 6 - 0.3);
			\node at (\a1 + 0.05, -\b / 6 + 0.5) {$S_{2,2}^{3,4}$};
			\draw [->] (\a2, -\b / 6 + 0.1) -- (\a1, -\b / 6 + 0.3);
			\node at (\a2 - 0.05, -\b / 6 + 0.5) {$S_{2,2}^{3,3}$};
			\draw [->] (\a2 - 0.2, -\b / 6 + 0.1) -- (\a2 - 0.1, -\b / 6 + 0.3);
			\node at (\a3 + 0.05, -\b / 6 - 0.55) {$S_{2,2}^{3,2}$};
			\draw [->] (\a3 + 0.05, -\b / 6 - 0.1) -- (\a3 + 0.05, -\b / 6 - 0.3);
			\node at (\a4, -\b / 6 - 0.55) {$S_{2,2}^{3,1}$};
			\draw [->] (\a4, -\b / 6 - 0.1) -- (\a4, -\b / 6 - 0.3);
			
		\end{tikzpicture}
		
		\caption{Level 3}\label{fig-segment-level3}
	\end{subfigure}
	
	\caption{Segments obtained from the separators. 
    (a): Separator $S_{1,1}^{1,1}$  splits the graph into subgraphs $G_{1,1}$ and $G_{1,2}$.
    (b) and (c): The graph is further separated into 4 subgraphs $G_{2,1}$, $G_{2,2}$, $G_{2,3}$, $G_{2,4}$. Meanwhile, separator $S^{1,1}_{1,1}$ is also split into three regular segments $S^{2,1}_{1,1}$, $S^{2,3}_{1,1}$, $S^{2,5}_{1,1}$ and two junction segments (in red) $S^{2,2}_{1,1}$, $S^{2,4}_{1,1}$,  by the new segments $S^{2,1}_{2,1}$ and $S^{2,1}_{2,2}$.
    (d): Segments $S^{2,1}_{2,1}$ and $S^{2,1}_{2,2}$ are further split into finer segments at higher level.
    }
    \label{fig-segment_generate}
\end{figure}

\subsection{Structure of segments}
\label{section-segment in ND}
In this section, we provide details on how to obtain segment information. Recall that in Section 2, we denote $S_{l,i}$ the $i$-th separator at level $l$. However, this notation is insufficient for constructing the tree structure of segments. To distinguish different segments, particularly their origins and how they divide, we introduce the notation $S_{l,i}^{j,k}$ to denote the $k$-th segment that belongs to the separator $S_{l,i}$, and is divided at level $j$.

As shown in \fig\ref{fig-segment_generate}, suppose that a graph $G=(V,E)$ is separated into two subgraphs $G_{1,1}$ and $G_{1,2}$ by separator $S_{1,1}$ at level 1, which is denoted as segment $S_{1,1}^{1,1}$ in our notation.
At level 2,  segment $S_{1,1}^{1,1}$ is further divided into $S_{1,1}^{2,1}\cup S_{1,1}^{2,2}\cup S_{1,1}^{2,3}$ by the separator $S_{2,1}$, or $S_{2,1}^{2,1}$ by segment notation.
Specifically, segment $S_{2,1}^{2,1}$ divides  $S_{1,1}^{1,1}$ into three segments: two disjoint \textit{regular} segments $S_{1,1}^{2,1}$ and $S_{1,1}^{2,3}$, and one \textit{junction} segment $S_{1,1}^{2,2}$.
In subgraph $G_{1,2}$, the separator $S_{2,2}$, also denoted as segment $S_{2,2}^{2,1}$, further splits the segment $S_{1,1}^{2,3}$ into three smaller segments: $S_{1,1}^{2,3}\cup S_{1,1}^{2,4}\cup S_{1,1}^{2,5}$.
Here, the same notation $S_{1,1}^{2,3}$ is used to denote one of the smaller segment. At level 3, segment $S_{2,i}^{2,1}$ is further divided into $\cup_{k=1}^5S_{2,i}^{3,k}$ with $i=1,2$, as shown in \fig\ref{fig-segment-level3}. Following this procedure, all segments can be recursively split into smaller ones, enabling the construction of a tree structure by indexing each segment appropriately.

In the end, the algorithm for partitioning the graph $G$ by separators with segment information is outlined in \algo\ref{alg-nested 2D}.

\section{Sparsification based on low rank approximations}\label{sec4}

In this section, we detail the transformation process mentioned in Section 2, with a specific focus on identifying the low rank structures using \textit{interpolative decomposition} (ID) \cite{ID-LowRank}. 
		
\subsection{Interpolative Decomposition}
\label{section-intro ID}
Given an $m \times n$ matrix $B$ with numerical rank $k$ under the tolerance $\eps$, i.e.,  $\sigma_k(B)\ge \eps>\sigma_{k+1}(B)$, where $\sigma_k(B)$ is the $k$-th singular value of $B$, one can apply the strong rank-reveal QR (RRQR)  factorization to $B$~\cite{Strong_RRQR}. It yields: 
\begin{equation}
    B \Pi = Q
    \begin{bmatrix}
        R_1 & R_2 \\
          & R_{\eps}
    \end{bmatrix}
    = Q 
    \begin{bmatrix}
        R_1 & R_2 \\
        & \O(\eps)
    \end{bmatrix},
\label{equ-RRQR 1}
\end{equation}
where $\Pi$ is a permutation matrix, which in most of the time is the same as the permutation generated by column-pivoted QR (CPQR)~\cite{MatrixComputation}. Here, $R_1$ is an upper triangular matrix whose last singular value is greater than $\O(\eps)$ and $Q$ is an orthogonal matrix of the form 
$\begin{bmatrix}
    Q_s & Q_r
\end{bmatrix}$, such that:
\begin{equation}
    B \Pi = Q_s \left[ R_1 \ R_2 \right] + \O(\eps)
    = Q_s R_1 \left[ I \quad R_1^{-1} R_2 \right] + \O(\eps)
    = B_s \left[ I \ T \right] + \O(\eps),
\label{equ-RRQR 2}
\end{equation}
where $T = R_1^{-1} R_2$, $\Pi = \left[ \Pi_s \quad \Pi_r \right]$, $\Pi_s$ represents the first $k$ columns of $\Pi$, and $B_s = Q_s R_1 = B \Pi_s$ is the \textit{skeleton} of $B$. The factorization in the form \eqref{equ-RRQR 2} is called the interpolative decomposition of $B$ \cite{ID-LowRank}, which satisfies: 
\begin{equation}
    B 
    \underbrace{
        \begin{bmatrix}
            \Pi_r & \Pi_s
        \end{bmatrix}
        \begin{bmatrix}
            I  & \\
            -T & I
        \end{bmatrix}
    }_{L_t}
    = 
    \begin{bmatrix}
        0 & B_s
    \end{bmatrix}+\O(\eps).
\label{equ-B}
\end{equation}

The transformation $U_t$ is generated through a similar procedure when the matrix $B$ is transposed. Compared to SVD, the interpolative decomposition provides a cost-effective alternative to approximate the low rank structure. Its worst case computational cost is $\O(mn^2)$, but typically is $\O(mnk)$, which is comparable to the Gram–Schmidt algorithm. Despite this, the computational cost remains significant and reducing it is highly desirable,  even if only the constant prefactor in the complexity can be reduced, as it dominates the cost of sparse decomposition.
		
\subsection{Accelerating the interpolative decomposition}
 Note that in the QR factorization~\eqref{equ-RRQR 2}, we only need the permutation matrix $\Pi$ and the matrix $T$ with $k$ ($k \ll m$) rows. 	Let $\Phi$ be an $h \times m$ matrix, referred to as the \textit{sampling} matrix. By defining $Y = \Phi B$ and applying RRQR to $Y$, we obtain:
\begin{equation}
    Y \Pi = Y \Pi_s \left[ I \ T \right] + \O(\eps),
\end{equation}
which is
\begin{equation}
    \Phi B \Pi = Y \Pi_s \left[ I \ T \right] + \O(\eps) = \Phi B_s \left[ I \ T \right] + \O(\eps).
\label{equ-RRQR sample}
\end{equation}
If $\Phi$ is invertible, in which case $h = m$, one can verify that 
\begin{equation}
    B \Pi = B_s \left[ I \ T \right] + \O(\eps),
\end{equation}
which has the same form as equation \eqref{equ-RRQR 2}. However, in practice, matrix $\Phi$ is typically of size $h < m$. To accelerate the interpolative decomposition, the challenge is to find an efficient method for selecting $\Phi$, which we call sampling matrix.
	
\subsubsection{Randomized Sampling}
\label{section-Random Sampling}
There are several methods available to determine the matrix $\Phi$. One approach is randomized sampling, which makes use of the properties of random matrices \cite{RandomMatrix}.  This method uses a Gaussian random matrix $G$ as the sampling matrix (i.e., $\Phi = G$) and then computes the $QR$ factorization of the resulting matrix $GB$ to efficiently identify the low rank structure of a dense matrix. That is: 
\begin{equation}
G B \Pi = Q R.
\end{equation}
The method is simple and elegant. However, it is purely algebraic and does not take into account the underlying geometric information of the matrix. As a result, it may not always be the most optimal choice.
	
\subsubsection{Sampling by FMM}
\label{section-Ideas from FMM}
The Fast Multipole Method (FMM) \cite{FMM-lecture, FMM-lecture2} is designed to accelerate matrix-vector multiplication in the $N$-body interactions by partitioning the interactions into near-field and far-field components. It hierarchically approximate the far-field interaction using low rank structures. FMM is particularly effective for linear elliptic PDEs because the kernels of various integral equations (such as Laplace or Helmholtz) are translation-invariant and decay with increasing distance.

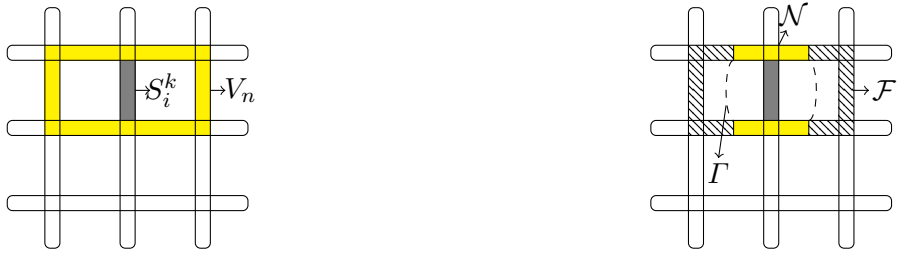
\begin{figure}[tb]
    \centering
    \begin{subfigure} [t] {0.48\textwidth}
        \centering
        \begin{tikzpicture}[scale=1, 
            pattern_slash/.style={pattern=north east lines},
            pattern_Rslash/.style={pattern=north west lines},]
            
            \fill[gray] (0.1,2.1) rectangle (-0.1,2.9);
            
            \fill[yellow] (-1.1,1.9) rectangle (1.1,2.1);
            \fill[yellow] (-1.1,2.9) rectangle (1.1,3.1);
            \fill[yellow] (-1.1,2.1) rectangle (-0.9,2.9);
            \fill[yellow] (1.1,2.1) rectangle (0.9,2.9);
            
            \draw [rounded corners=2] (-0.1,0.4) rectangle (0.1,3.6);
            
            \draw [rounded corners=2] (-1.6,0.9) rectangle (1.6,1.1);
            \draw [rounded corners=2] (-1.6,1.9) rectangle (1.6,2.1);
            \draw [rounded corners=2] (-1.6,2.9) rectangle (1.6,3.1);
            
            \draw [rounded corners=2] (-0.9,0.4) rectangle (-1.1,3.6);
            \draw [rounded corners=2] (0.9,0.4) rectangle (1.1,3.6);
            
            \node at (0.45, 2.5) {$S^k_i$};
            \draw [->] (0.1, 2.5) -- (0.3, 2.5);
            \node at (1.5, 2.5) {$V_n$};
            \draw [->] (1.1, 2.5) -- (1.3, 2.5);
            
        \end{tikzpicture}
        
        \caption{Edge relation between segment $S^k_i$ and its neighbor $V_n$}\label{fig-FMM proxy curve_a}
    \end{subfigure}
    \quad
    \begin{subfigure} [t] {0.48\textwidth}
        \centering
        \begin{tikzpicture}[scale=1, 
            pattern_slash/.style={pattern=north east lines},
            pattern_Rslash/.style={pattern=north west lines},]
            
            \fill[gray] (0.1,2.1) rectangle (-0.1,2.9);
            
            \fill[yellow] (-0.5,1.9) rectangle (0.5,2.1);
            \fill[yellow] (-0.5,2.9) rectangle (0.5,3.1);
            \fill [pattern_Rslash] (-1.1, 1.9) rectangle (-0.5, 2.1);
            \fill [pattern_Rslash] (-1.1, 2.9) rectangle (-0.5, 3.1);
            \fill [pattern_Rslash] (-1.1, 2.1) rectangle (-0.9, 2.9);
            \fill [pattern_Rslash] (1.1, 1.9) rectangle (0.5, 2.1);
            \fill [pattern_Rslash] (1.1, 2.9) rectangle (0.5, 3.1);
            \fill [pattern_Rslash] (1.1, 2.1) rectangle (0.9, 2.9);
            
            \draw (-0.5,1.9) -- (-0.5,2.1);
            \draw (0.5,1.9) -- (0.5,2.1);
            \draw (-0.5,2.9) -- (-0.5,3.1);
            \draw (0.5,2.9) -- (0.5,3.1);
            
            \draw [rounded corners=2] (-0.1,0.4) rectangle (0.1,3.6);
            
            \draw [rounded corners=2] (-1.6,0.9) rectangle (1.6,1.1);
            \draw [rounded corners=2] (-1.6,1.9) rectangle (1.6,2.1);
            \draw [rounded corners=2] (-1.6,2.9) rectangle (1.6,3.1);
            
            \draw [rounded corners=2] (-0.9,0.4) rectangle (-1.1,3.6);
            \draw [rounded corners=2] (0.9,0.4) rectangle (1.1,3.6);
            
            \draw[dashed] (0.5,2.1) arc (-90:90:0.1 and 0.4);
            \draw[dashed] (-0.5,2.9) arc (90:270:0.1 and 0.4);
            
            \node at (0.3, 3.5) {$\mathcal{N}$};
            \draw [->] (0.1, 3.1) -- (0.2, 3.3);
            \node at (1.5, 2.5) {$\mathcal{F}$};
            \draw [->] (1.1, 2.5) -- (1.3, 2.5);
            \node at (-0.7
            , 1.4) {$\Gamma$};
            \draw [->] (-0.6, 2.3) -- (-0.7, 1.6);
            
        \end{tikzpicture}
        
        \caption{Far-field $\mathcal{F}$ is replaced by a proxy curve $\Gamma$}\label{fig-FMM proxy curve_b}
    \end{subfigure}
    
    \caption{An FMM approach to handling the far-field $\mathcal{F}$. 
    (a): The edge relation between segment $S^k_i$ and its neighboring segments $V_n$. 
    (b): The neighbor of $S^k_i$ is partitioned as $V_n = \mathcal{N} \cup \mathcal{F}$ where $\mathcal{N}$ represents the near-field and $\mathcal{F}$ can be replaced by a proxy curve $\Gamma$.}
    \label{fig-FMM proxy curve}
\end{figure}

In terms of numerical algebra, consider a sub-matrix $A_{n,i}$ of $A$, where $A_{n,i}$ represents the edge relations between segment $S^k_i$ and its neighboring set $V_n$, as shown in \fig\ref{fig-FMM proxy curve_a}. This sub-matrix $A_{n,i}$ can be compressed via low rank approximation with tolerance $\eps$, such that,
\begin{equation}
    P^\top A_{n,i}
    = \begin{bmatrix}
        A_{\mathcal{N},i} \\ A_{\mathcal{F},i}
    \end{bmatrix}
    = \begin{bmatrix}
        I & \\
          & M_{\mathcal{F},\Gamma}
    \end{bmatrix}
    \begin{bmatrix}
        A_{\mathcal{N},i} \\ Y_{\Gamma,i}
    \end{bmatrix} + \O(\eps)
\label{equ-proxy curve}
\end{equation}
where $P^\top$ is a permutation matrix that partitions $V_n$ into near-field ($\mathcal{N}$) and far-field ($\mathcal{F}$). 
As illustrated in \fig\ref{fig-FMM proxy curve_b}, $\Gamma$ represents the set of indices corresponding to the dashed curve  with  size $\abs(\Gamma) < \abs(\mathcal{F})$.
The edge relations between vertex set $i$ and the far field $\mathcal{F}$, represented by $A_{\mathcal{F},i}$, are approximated by the relation between $i$ and the virtual proxy vertex set $\Gamma$ represented by $Y_{\Gamma,i}$ \cite{FMM-skeleton}.
Therefore, we can take $\Gamma$ as the sampling matrix.
However, the size of $\Gamma$ highly depends on the distance between the set $i$ and its far-field vertex set $\mathcal{F}$, which increases rapidly if $\mathcal{F}$ is closely connected to $i$. 
		
\subsubsection{A hybrid sampling}
In this part, we combine the benefits of randomized sampling and the sampling idea from FMM. Consider a sub-matrix $A_{n,i}$ with numerical rank $k$ and let $P^\top$ be a permutation matrix in the form of (\ref{equ-proxy curve}).  Define the sampling matrix $\Phi$ as
\begin{equation}
    \Phi
    = \begin{bmatrix}
        I & \\
          & G
    \end{bmatrix}
    P^\top,
\label{equ-Phi}
\end{equation}
where $G$ is a Gaussian random matrix of size $h \times \abs(\mathcal{F})$.
Based on \cite{RandomMatrix}, for $A_{\mathcal{F},i}$ in \eqref{equ-proxy curve}, there exists an orthonormal matrix $Q_\mathcal{F} \in \real^{\abs(\mathcal{F}), h}$ with $h = k + p$ ($p=5$ is sufficient in practice) such that  
\begin{equation}
 A_{\mathcal{F},i} - A_{\mathcal{F},i} Q_\mathcal{F} Q_\mathcal{F}^\top   = \O(\sigma_{k+1}), 
\end{equation}
where $\sigma_k$ is the $k$-th singular value of $A_{\mathcal{F},i}$ and the QR factorization yields $(G A_{\mathcal{F},i})^\top = Q_\mathcal{F} R_\mathcal{F}$. 
In other words, $A_{\mathcal{F},i} \approx A_{\mathcal{F},i} Q_\mathcal{F} Q_\mathcal{F}^\top = \alpha_\mathcal{F} Q_\mathcal{F}^\top$, where $\alpha_\mathcal{F} = A_{\mathcal{F},i} Q_\mathcal{F}$.
Therefore, 
\begin{equation}
\begin{bmatrix}
    A_{\mathcal{N},i} \\ A_{\mathcal{F},i}
\end{bmatrix} = 
\begin{bmatrix}
    I & \\ & \alpha_\mathcal{F} R_\mathcal{F}^{-T}
\end{bmatrix}
\begin{bmatrix}
    A_{\mathcal{N},i} \\ R_\mathcal{F}^\top Q_\mathcal{F}^\top
\end{bmatrix} + \O(\sigma_{k+1}) = 
\begin{bmatrix}
    I & \\ & \alpha_\mathcal{F} R_\mathcal{F}^{-T}
\end{bmatrix}
\begin{bmatrix}
    A_{\mathcal{N},i} \\ G A_{\mathcal{F},i}
\end{bmatrix} + \O(\sigma_{k+1}).
\end{equation}
As a result, the interpolative decomposition can be applied to the matrix sampled by the sampling matrix (\ref{equ-Phi}), denoted by $Y_{n,i} = \Phi A_{n,i}$, rather than directly to $A_{n,i}$.

Compared to just using randomized sampling given in subsection \ref{section-Random Sampling}, this approach requires a smaller Gaussian random matrix. Additionally, there is no need to explicitly compute the proxy mapping $M_{\mathcal{F},\Gamma}$ as required in the FMM, since random matrix $G$ effectively replaces the proxy mapping in revealing the rank of interpolative decomposition.

\subsection{SpaLU}
The complete spaLU algorithm to factorize the matrix $A$ in \equ(\ref{equ-Ax=b}) consists of the following steps: elimination of interiors, sparsification, elimination, and merging of the segments.

First, we eliminate interiors for all the subgraphs  at the leaf level via straightforward LU factorization, leaving a remainder graph that consists solely of separators. This step yields
\begin{equation}
    L_{(L+1)}^{-1} A U_{(L+1)}^{-1} = \begin{bmatrix}
        I & \\ & A^{(L)}
    \end{bmatrix},
\end{equation}
where $A \in \complex^{N,N}, \ A^{(L)} \in \complex^{N_L,N_L}$, and $N_L$ is the number of vertices on the separators. Note that the first $N-N_L$ rows and columns of $A$ form a block diagonal matrix and can be factorized in parallel.

\begin{figure}[tb]
    \centering
    \begin{subfigure} [t] {0.3\textwidth}
        \centering
        \begin{tikzpicture}[scale=0.4, 
            pattern_slash/.style={pattern=north east lines} ]
            
            \foreach \i in {0,2,...,10}
            { \draw (0, -\i) -- (10, -\i); }
            \foreach \i in {0,2,...,10}
            { \draw (\i, 0) -- (\i, -10); }
            
            \filldraw [fill=white] (0,0) rectangle (4,-4);
            \foreach \i in {1,...,10}
            { \filldraw[fill=gray] (\i-1,-\i+1) rectangle (\i,-\i); }
            
            \tikzmath{ \i = 1; }; \foreach \j in {5,6,10}
            { \fill[pattern_slash] (\j-1, -\i+1) rectangle (\j, -\i); 
              \fill[pattern_slash] (\i-1, -\j+1) rectangle (\i, -\j); }
            \draw (4,-1) -- (6,-1); \draw (1,-4) -- (1,-6);
            \draw (8,-1) -- (10,-1); \draw (1,-8) -- (1,-10);
            \draw (9,0) -- (9,-4); \draw (0,-9) -- (4,-9);
            \tikzmath{ \i = 2; }; \foreach \j in {5,6,9}
            { \fill[pattern_slash] (\j-1, -\i+1) rectangle (\j, -\i); 
              \fill[pattern_slash] (\i-1, -\j+1) rectangle (\i, -\j); }
            \tikzmath{ \i = 3; }; \foreach \j in {7,8,10}
            { \fill[pattern_slash] (\j-1, -\i+1) rectangle (\j, -\i); 
              \fill[pattern_slash] (\i-1, -\j+1) rectangle (\i, -\j); }
            \draw (6,-3) -- (10,-3); \draw (3,-6) -- (3,-10);
            \tikzmath{ \i = 4; }; \foreach \j in {7,8,9}
            { \fill[pattern_slash] (\j-1, -\i+1) rectangle (\j, -\i); 
                \fill[pattern_slash] (\i-1, -\j+1) rectangle (\i, -\j); }
            \tikzmath{ \i = 6; }; \foreach \j in {9,10}
            { \fill[pattern_slash] (\j-1, -\i+1) rectangle (\j, -\i); 
                \fill[pattern_slash] (\i-1, -\j+1) rectangle (\i, -\j); }
            \draw (8,-5) -- (10,-5); \draw (5,-8) -- (5,-10);
            \tikzmath{ \i = 7; }; \foreach \j in {9,10}
            { \fill[pattern_slash] (\j-1, -\i+1) rectangle (\j, -\i); 
                \fill[pattern_slash] (\i-1, -\j+1) rectangle (\i, -\j); }
            \draw (8,-7) -- (10,-7); \draw (7,-8) -- (7,-10);
        \end{tikzpicture}
        
        \caption{Schur complement of $A^{(l)}$}\label{fig-sp-lu (a)}
    \end{subfigure}
    \
    \begin{subfigure} [t] {0.3\textwidth}
        \centering
        \begin{tikzpicture}[scale=0.4, 
            pattern_Rslash/.style={pattern=north west lines} ]
            
            \foreach \i in {0,2,...,10}
            { \draw (0, -\i) -- (10, -\i); }
            \foreach \i in {0,2,...,10}
            { \draw (\i, 0) -- (\i, -10); }
            
            \filldraw [fill=white] (0,0) rectangle (4,-4);
            \foreach \i in {1,...,10}
            { \filldraw[fill=gray] (\i-1,-\i+1) rectangle (\i,-\i); }
            
            \tikzmath{ \i = 1; }; \foreach \j in {5,6,10}
            { \filldraw[pattern_Rslash] (\j-0.5, -\i+0.5) rectangle (\j, -\i); 
                \filldraw[pattern_Rslash] (\i-0.5, -\j+0.5) rectangle (\i, -\j); }
            \draw (4,-1) -- (6,-1); \draw (1,-4) -- (1,-6);
            \draw (8,-1) -- (10,-1); \draw (1,-8) -- (1,-10);
            \draw (9,0) -- (9,-4); \draw (0,-9) -- (4,-9);
            \tikzmath{ \i = 2; }; \foreach \j in {5,6,9}
            { \filldraw[pattern_Rslash] (\j-0.5, -\i+0.5) rectangle (\j, -\i); 
                \filldraw[pattern_Rslash] (\i-0.5, -\j+0.5) rectangle (\i, -\j); }
            \tikzmath{ \i = 3; }; \foreach \j in {7,8,10}
            { \filldraw[pattern_Rslash] (\j-0.5, -\i+0.5) rectangle (\j, -\i); 
                \filldraw[pattern_Rslash] (\i-0.5, -\j+0.5) rectangle (\i, -\j); }
            \draw (6,-3) -- (10,-3); \draw (3,-6) -- (3,-10);
            \tikzmath{ \i = 4; }; \foreach \j in {7,8,9}
            { \filldraw[pattern_Rslash] (\j-0.5, -\i+0.5) rectangle (\j, -\i); 
                \filldraw[pattern_Rslash] (\i-0.5, -\j+0.5) rectangle (\i, -\j); }
            \tikzmath{ \i = 6; }; \foreach \j in {9,10}
            { \filldraw[pattern_Rslash] (\j-0.5, -\i+0.5) rectangle (\j, -\i); 
                \filldraw[pattern_Rslash] (\i-0.5, -\j+0.5) rectangle (\i, -\j); }
            \draw (8,-5) -- (10,-5); \draw (5,-8) -- (5,-10);
            \tikzmath{ \i = 7; }; \foreach \j in {9,10}
            { \filldraw[pattern_Rslash] (\j-0.5, -\i+0.5) rectangle (\j, -\i); 
                \filldraw[pattern_Rslash] (\i-0.5, -\j+0.5) rectangle (\i, -\j); }
            \draw (8,-7) -- (10,-7); \draw (7,-8) -- (7,-10);
        \end{tikzpicture}
        
        \caption{After sparsification}\label{fig-sp-lu (b)}
    \end{subfigure}
    \
    \begin{subfigure} [t] {0.3\textwidth}
        \centering
        \begin{tikzpicture}[scale=0.4, 
            pattern_Rslash/.style={pattern=north west lines} ]
            
            \foreach \i in {0,2,...,10}
            { \draw (0, -\i) -- (10, -\i); }
            \foreach \i in {0,2,...,10}
            { \draw (\i, 0) -- (\i, -10); }
            
            \filldraw [fill=white] (0,0) rectangle (4,-4);
            \foreach \i in {1,...,10}
            { \draw (\i-1, -\i+1) rectangle (\i, -\i); }
            \draw[thick] (0,0) -- (10,-10);
            
            \tikzmath{ \i = 1; }; \foreach \j in {5,6,10}
            { \filldraw[pattern_Rslash] (\j-0.5, -\i+0.5) rectangle (\j, -\i); 
                \filldraw[pattern_Rslash] (\i-0.5, -\j+0.5) rectangle (\i, -\j); }
            \draw (4,-1) -- (6,-1); \draw (1,-4) -- (1,-6);
            \draw (8,-1) -- (10,-1); \draw (1,-8) -- (1,-10);
            \draw (9,0) -- (9,-4); \draw (0,-9) -- (4,-9);
            \tikzmath{ \i = 2; }; \foreach \j in {5,6,9}
            { \filldraw[pattern_Rslash] (\j-0.5, -\i+0.5) rectangle (\j, -\i); 
                \filldraw[pattern_Rslash] (\i-0.5, -\j+0.5) rectangle (\i, -\j); }
            \tikzmath{ \i = 3; }; \foreach \j in {7,8,10}
            { \filldraw[pattern_Rslash] (\j-0.5, -\i+0.5) rectangle (\j, -\i); 
                \filldraw[pattern_Rslash] (\i-0.5, -\j+0.5) rectangle (\i, -\j); }
            \draw (6,-3) -- (10,-3); \draw (3,-6) -- (3,-10);
            \tikzmath{ \i = 4; }; \foreach \j in {7,8,9}
            { \filldraw[pattern_Rslash] (\j-0.5, -\i+0.5) rectangle (\j, -\i); 
                \filldraw[pattern_Rslash] (\i-0.5, -\j+0.5) rectangle (\i, -\j); }
            \tikzmath{ \i = 6; }; \foreach \j in {9,10}
            { \filldraw[pattern_Rslash] (\j-0.5, -\i+0.5) rectangle (\j, -\i); 
                \filldraw[pattern_Rslash] (\i-0.5, -\j+0.5) rectangle (\i, -\j); }
            \draw (8,-5) -- (10,-5); \draw (5,-8) -- (5,-10);
            \tikzmath{ \i = 7; }; \foreach \j in {9,10}
            { \filldraw[pattern_Rslash] (\j-0.5, -\i+0.5) rectangle (\j, -\i); 
                \filldraw[pattern_Rslash] (\i-0.5, -\j+0.5) rectangle (\i, -\j); }
            \draw (8,-7) -- (10,-7); \draw (7,-8) -- (7,-10);
        \end{tikzpicture}

        \caption{After LU factorization}
       \label{fig-sp-LU (c)}
    \end{subfigure}
    \
    \begin{subfigure} [t] {0.3\textwidth}
        \centering
        \begin{tikzpicture}[scale=0.4, 
            pattern_slash/.style={pattern=north east lines} ]
            
            \foreach \i in {0,2,...,10}
            { \draw (0, -\i) -- (10, -\i); }
            \foreach \i in {0,2,...,10}
            { \draw (\i, 0) -- (\i, -10); }
            
            \filldraw [fill=white] (0,0) rectangle (4,-4);
            \draw[thick] (0,0) -- (10,-10);
            
            \tikzmath{ \i = 5; }; \foreach \j in {6,9,10}
            { \filldraw[pattern_slash] (\j-0.5, -\i+0.5) rectangle (\j, -\i); 
                \filldraw[pattern_slash] (\i-0.5, -\j+0.5) rectangle (\i, -\j); }
            \tikzmath{ \i = 6; }; \foreach \j in {5,9,10}
            { \filldraw[pattern_slash] (\j-0.5, -\i+0.5) rectangle (\j, -\i); 
                \filldraw[pattern_slash] (\i-0.5, -\j+0.5) rectangle (\i, -\j); }
            \draw (8,-5) -- (10,-5); \draw (5,-8) -- (5,-10);
            \tikzmath{ \i = 7; }; \foreach \j in {8,9,10}
            { \filldraw[pattern_slash] (\j-0.5, -\i+0.5) rectangle (\j, -\i); 
                \filldraw[pattern_slash] (\i-0.5, -\j+0.5) rectangle (\i, -\j); }
            \draw (8,-7) -- (10,-7); \draw (7,-8) -- (7,-10);
            \tikzmath{ \i = 8; }; \foreach \j in {7,9,10}
            { \filldraw[pattern_slash] (\j-0.5, -\i+0.5) rectangle (\j, -\i); 
                \filldraw[pattern_slash] (\i-0.5, -\j+0.5) rectangle (\i, -\j); }
            \foreach \i in {5,...,10}
            { \filldraw[fill=gray] (\i-0.5,-\i+0.5) rectangle (\i,-\i); }
        \end{tikzpicture}

        \caption{After eliminating separators at level $l$}
        \label{fig-sp-LU (d)}
    \end{subfigure}
    \
    \begin{subfigure} [t] {0.3\textwidth}
        \centering
        \begin{tikzpicture}[scale=0.4, 
            pattern_slash/.style={pattern=north east lines} ]
            
            \draw (-7,7) rectangle (3,-3);
            \draw[thick] (-7,7) -- (3, -3);
            \draw (0,0) rectangle (3,-3);
            
            \tikzmath{ \i = 1; }; \foreach \j in {3}
            { \filldraw[pattern_slash] (\j-1, -\i+1) rectangle (\j, -\i); 
                \filldraw[pattern_slash] (\i-1, -\j+1) rectangle (\i, -\j); }
            \tikzmath{ \i = 2; }; \foreach \j in {3}
            { \filldraw[pattern_slash] (\j-1, -\i+1) rectangle (\j, -\i); 
                \filldraw[pattern_slash] (\i-1, -\j+1) rectangle (\i, -\j); }
            \foreach \i in {1,2}
            { \filldraw[fill=gray] (\i-1,-\i+1) rectangle (\i,-\i); }
            \foreach \i in {5,6}
            { \filldraw[fill=gray] (0.5*\i-0.5,-0.5*\i+0.5) rectangle (0.5*\i,-0.5*\i); }
            
            \draw (2.5, 0) -- (2.5,-2); \draw (0, -2.5) -- (2, -2.5);
            
        \end{tikzpicture}
        
        \caption{After merging the segments}
        \label{fig-sp-LU (e)}
    \end{subfigure}
    
    \caption{An illustration of spaLU in terms of matrix pattern, which begins with the elimination of interior vertices at the leaf level, followed by sparsification, elimination, and merging of segments. Here regions in gray, slashed or back-slashed represent dense sub-matrices and areas in white refer to zero-matrices. Diagonal with thick line denotes identity matrix.}
    \label{fig-sp-LU matrix}
\end{figure}
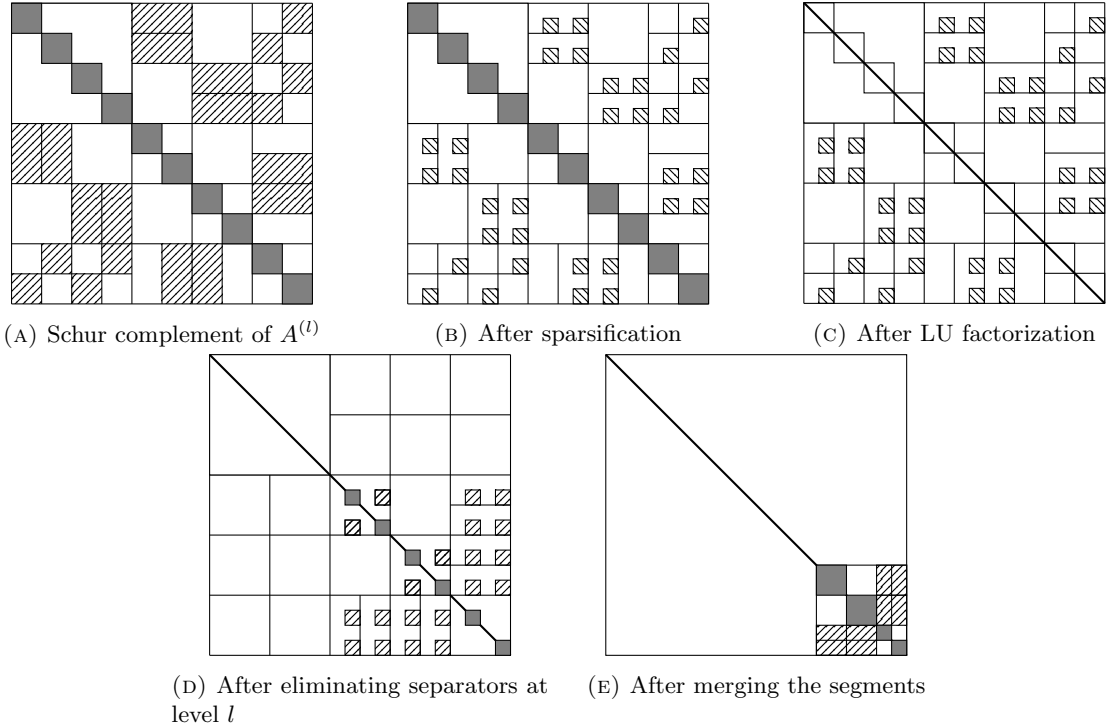

Next, we focus on the remaining graph corresponding to the Schur complement $A^{(L)}$,  as shown in \fig\ref{fig-sp-lu (a)}, which is associated with the separator graph described in subsection \ref{section_segments}. Since it is factorized in a recursive way, let us assume the separators at level $l+1$ has been eliminated and the remaining matrix becomes $A^{(l)}$. We then sparsify each segment $S_{i,j}^{l,k}$ in the rest of the graph via interpolative decomposition, as described in subsection \ref{section-intro Sparsification} and illustrated in \fig\ref{fig-sp-lu (b)}, which yields:
\begin{equation}
    A_{sp}^{(l)} = U_T^{(l)} A^{(l)} L_T^{(l)} + \O(\eps).
\label{equ-U_T * A * LT}
\end{equation}
Apply LU factorization to the submatrices associated with the sparsified segments as shown in \fig\ref{fig-sp-LU (c)}, and use $Z_l^{(l)}$ and $Z_u^{(l)}$ to eliminate the segments $S_{l,j}^{i,k}$, as shown in \fig\ref{fig-sp-LU (d)}, which yields
\begin{equation}
Z_l^{(l)} \begin{bmatrix}
L_{l,1}^{-1} & & \\
& L_{l,2}^{-1} & \\
& & \ddots
\end{bmatrix}
A_{sp}^{(l)}
\begin{bmatrix}
U_{l,1}^{-1} & & \\
& U_{l,2}^{-1} & \\
& & \ddots
\end{bmatrix} Z_u^{(l)} = 
L_{(l)}^{-1} A_{sp}^{(l)} U_{(l)}^{-1} = A_{lu}^{(l)}.
\label{equ-L^-1 * A * U^-1}
\end{equation}
	
Then, we permute the matrix to merge the indices, corresponding to the skeletons after sparsification, into the last of rows and columns, as illustrated in \fig\ref{fig-sp-LU (e)}. This results in: 
\begin{equation}
    A^{(l-1)} = E_{(l)}^\top A_{lu}^{(l)} E_{(l)}
\end{equation}
where $E_{(l)}$ is a permute matrix.

In the end, $A^{(L)}$ is factorized into
\begin{equation}
    W A^{(L)} V \approx I,
    \end{equation}
with
\begin{equation}    
    W = \prod_{l=1}^{L}\left( E^\top_{(l)} L_{(l)}^{-1} U_T^{(l)} \right),\
    V = \prod_{l=L}^{1}\left( L_T^{(l)} U_{(l)}^{-1} E_{(l)} \right).
\label{equ-WV}
\end{equation}
Therefore, it holds
\begin{equation}
\begin{aligned}
    A &\approx \tilde{W}^{-1} \tilde{V}^{-1}, \\
    \tilde{W} &= 
    \begin{bmatrix}
        I & \\ & W
    \end{bmatrix}
    L_{(L+1)}^{-1}, \\
    \tilde{V} &= 
    U_{(L+1)}^{-1} 
    \begin{bmatrix}
        I & \\ & V
    \end{bmatrix}.
\end{aligned}
\end{equation}
In the end, we have $\mathcal{A}  = \mathfrak{I}A \mathfrak{I}^\top $. The whole process is summarized in \algo\ref{alg-sp-LU}.

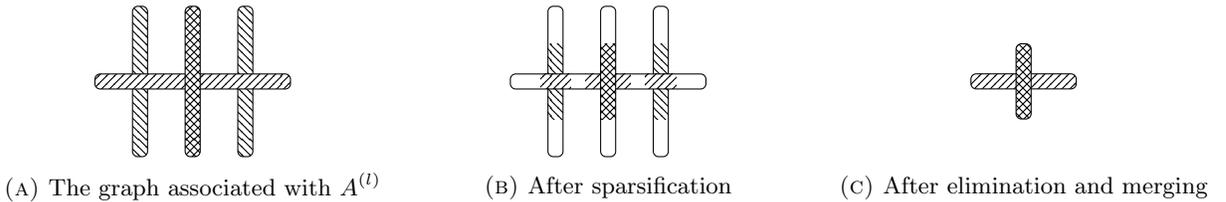
\begin{figure}[tb]
    \centering
    \begin{subfigure} [t] {0.3\textwidth}
        \centering
        \begin{tikzpicture}[scale=1, 
            pattern_slash/.style={pattern=north east lines},
            pattern_Rslash/.style={pattern=north west lines} ]
            
            \filldraw[pattern_slash] [rounded corners=2] (-0.1,1) rectangle (0.1,-1);
            \fill[pattern_Rslash] [rounded corners=2] (-0.1,1) rectangle (0.1,-1);
            
            \filldraw[pattern_slash] [rounded corners=2] (-0.1,-0.1) -- (-1.3,-0.1) -- (-1.3,0.1) -- (-0.1,0.1);
            \filldraw[pattern_slash] [rounded corners=2] (0.1,-0.1) -- (1.3,-0.1) -- (1.3,0.1) -- (0.1,0.1);
            
            \filldraw[pattern_Rslash] [rounded corners=2] (-0.6,-0.1) -- (-0.6,-1) -- (-0.8,-1) -- (-0.8,-0.1);
            \filldraw[pattern_Rslash] [rounded corners=2] (0.6,-0.1) -- (0.6,-1) -- (0.8,-1) -- (0.8,-0.1);
            \filldraw[pattern_Rslash] [rounded corners=2] (-0.6,0.1) -- (-0.6,1) -- (-0.8,1) -- (-0.8,0.1);
            \filldraw[pattern_Rslash] [rounded corners=2] (0.6,0.1) -- (0.6,1) -- (0.8,1) -- (0.8,0.1);
            
        \end{tikzpicture}
        
        \caption{The graph associated with $A^{(l)}$}
    \end{subfigure}
    \quad
    \begin{subfigure} [t] {0.3\textwidth}
        \centering
        \begin{tikzpicture}[scale=1, 
            pattern_slash/.style={pattern=north east lines},
            pattern_Rslash/.style={pattern=north west lines} ]
            
            \draw [rounded corners=2] (-0.1,1) rectangle (0.1,-1);
            \fill [pattern_slash] (-0.1,0.5) rectangle (0.1,-0.5);
            \fill [pattern_Rslash] (-0.1,0.5) rectangle (0.1,-0.5);
            
            \draw [rounded corners=2] (-0.1,-0.1) -- (-1.3,-0.1) -- (-1.3,0.1) -- (-0.1,0.1);
            \fill [pattern_slash] (-0.1,-0.1) rectangle (-0.3,0.1);
            \fill [pattern_slash] (-0.5,-0.1) rectangle (-0.9,0.1);
            \draw [rounded corners=2] (0.1,-0.1) -- (1.3,-0.1) -- (1.3,0.1) -- (0.1,0.1);
            \fill [pattern_slash] (0.1,-0.1) rectangle (0.3,0.1);
            \fill [pattern_slash] (0.5,-0.1) rectangle (0.9,0.1);
            
            \draw [rounded corners=2] (-0.6,-0.1) -- (-0.6,-1) -- (-0.8,-1) -- (-0.8,-0.1);
            \fill [pattern_Rslash] (-0.6,-0.1) rectangle (-0.8,-0.5);
            \draw [rounded corners=2] (0.6,-0.1) -- (0.6,-1) -- (0.8,-1) -- (0.8,-0.1);
            \fill [pattern_Rslash] (0.6,-0.1) rectangle (0.8,-0.5);
            \draw [rounded corners=2] (-0.6,0.1) -- (-0.6,1) -- (-0.8,1) -- (-0.8,0.1);
            \fill [pattern_Rslash] (-0.6,0.1) rectangle (-0.8,0.5);
            \draw [rounded corners=2] (0.6,0.1) -- (0.6,1) -- (0.8,1) -- (0.8,0.1);
            \fill [pattern_Rslash] (0.6,0.1) rectangle (0.8,0.5);
            
        \end{tikzpicture}
        
        \caption{After sparsification}
    \end{subfigure}
    \quad
    \begin{subfigure} [t] {0.3\textwidth}
        \centering
        \begin{tikzpicture}[scale=1, 
            pattern_slash/.style={pattern=north east lines},
            pattern_Rslash/.style={pattern=north west lines}]
            
            \filldraw[pattern_slash] [rounded corners=2] (-0.1,0.5) rectangle (0.1,-0.5);
            \fill[pattern_Rslash] [rounded corners=2] (-0.1,0.5) rectangle (0.1,-0.5);
            
            \filldraw[pattern_slash] [rounded corners=2] (-0.1,-0.1) -- (-0.7,-0.1) -- (-0.7,0.1) -- (-0.1,0.1);
            \filldraw[pattern_slash] [rounded corners=2] (0.1,-0.1) -- (0.7,-0.1) -- (0.7,0.1) -- (0.1,0.1);
            
            \node at (0,-0.88) {};
        \end{tikzpicture}
        
        \caption{After elimination and merging}
    \end{subfigure}
    
    \caption{An illustration of spaLU in terms of graph, including the sparsification, elimination, and merging of segments. Here regions in backslash line represent segments of level $l$, slashed areas corresponds to level $l-1$, and cross region is associated with level $l-2$.}
    \label{fig-sp-LU graph}
\end{figure}

\begin{algorithm}[htb]
    \caption{SpaLU factorization algorithm}
    \label{alg-sp-LU}
    \begin{algorithmic}[1]
        \Require Nested ordered matrix $A$ with segment information $\mathcal{N}$, and the tolerance $\eps$.
        \Ensure Factorization term $L_{(L+1)}, U_{(L+1)}, W, V$.
        
        \State $A^{(L)}, L_{(L+1)}, U_{(L+1)} \gets $ the matrix after eliminating interiors of $A$.
        
        \For{level $l = L, \dots, 1$}
            \For{each segment $S_{i,j}^{l,k}$} \label{alg-sp-LU use S_rec}
                \State $A^{(l)} \gets$ sparsify $S_{i,j}^{l,k}$.
            \EndFor
            \State Get $U_T^{(l)}$ and $L_T^{(l)}$ through \equ(\ref{equ-U_T * A * LT}).
            
            \For{each segment $S_{i,j}^{l,k}$}
                \State $A^{(l)} \gets$ apply LU factorization to $S_{i,j}^{l,k}$.
            \EndFor
            \For{each segment $S_{l,j}^{i,k}$}
                \State $A^{(l)}\gets$ eliminate $S_{l,j}^{i,k}$.
            \EndFor
            \State Get $L_{(l)}$, and $U_{(l)}$ through \equ(\ref{equ-L^-1 * A * U^-1}).
            
            \State Permute $A^{(l)}$ by merging the sparsified segments, and get $E_{(l)}$.
        \EndFor
        
        \State \Return $L_{(L+1)}, U_{(L+1)}, W, V$ where $W, V$ are defined in \equ(\ref{equ-WV}).
    \end{algorithmic}
\end{algorithm}
	
\subsection{Stability Analysis}
In this section, we demonstrate that the spaLU algorithm is stable under certain conditions, both for low rank decomposition and symmetric decomposition.

\subsubsection{Stability of low rank decomposition}
Denote $E$ the perturbation error for $A$ after the spaLU decomposition. From equation \eqref{equ-spLU_lite}, it is clear that the error introduced by the low rank approximation is:
\begin{align}
	\norm(E) &= \norm(A-\mathcal{U}_L \mathcal{U}_{L-1}\cdots \mathcal{U}_1 \mathcal{L}_1 \cdots \mathcal{L}_{L-1} \mathcal{L}_L) \notag \\
	&\le  \prod_{l=1}^L(1 + \norm(T_{U_l})) (1 + \norm(T_{L_l})) \cdot \eps \norm(A),
\end{align}
where $T_{U_l}$ and $T_{L_l}$ are matrices of the form  $T$ given in equation (\ref{equ-RRQR 2}), associated with $\mathcal{U}_l$ and $\mathcal{L}_l$, respectively. The norm of $T$ in \equ(\ref{equ-RRQR 2}) may  become uncontrollable if $R_1$ is highly ill-conditioned. However, Gu and Eisenstat \cite{Strong_RRQR, ID-LowRank} provided a stable algorithm for the RRQR factorization (\ref{equ-RRQR 1}), which states: 
\begin{equation}
    \begin{cases}
        \sigma_k(R_1) &\geq \frac{1}{\sqrt{1+nk(n-k)}} \sigma_k(B) , \\
        \sigma_k(R_{\eps}) &\le \sqrt{1+nk(n-k)} \sigma_{k+1}(B), \\
        \norm(T)_F &\le \sqrt{nk(n-k)}.
    \end{cases}
    \label{equ-norm(T)}
\end{equation}
It implies that $\norm(T_{U_l}) =\O(1)$ and $\norm(T_{L_l})=\O(1)$, ensuring stability of the low rank decomposition.

\subsubsection{Stability of symmetric decomposition}
Consider the sub-matrix of $A$ in the form 
\begin{equation}
	A_{\mathrm{sub}} = \begin{bmatrix}
		A_{i,i} & A_{i,n} \\
		A_{n,i} & A_{n,n}
	\end{bmatrix}.
\end{equation}
If $A \in \complex^{N \times N}$ is symmetric, $A_{i,n}$ and $A_{i,n}^*$ share the same interpolative decomposition, allowing the use of $L D L^\top$ instead of LU factorization. 
Thus, we have  $U_{(l)} = L_{(l)}^*$ for $l = 0, \dots, L$, which implies that spaLU preserves symmetry, i.e., $\tilde{W} = \tilde{V}^*$.

Let  $E = A_{\mathrm{sub}} - A_{\mathrm{sub,sp}}$ be the perturbation, where $A_{\mathrm{sub,sp}}$ is the sparsified matrix. By the Weyl’s inequality, for any symmetric matrices $A$ and $B$, it holds
$\lambda_k(A) + \lambda_{\min}(B) \le \lambda_k(A+B) \le \lambda_k(A) + \lambda_{\max}(B)$, where $\lambda_k(\cdot)$ is the $k$-th eigenvalue. Therefore, 
$\abs({ \lambda_k(A+B) - \lambda_k(A) })  \le \norm(B)$. As a result, if $\norm(E) < \lambda_{\min}(A_{\mathrm{sub}})$, then $A_{\mathrm{sub,sp}}$ remains symmetric and stable.		
	
\subsection{A note on factorizing unsymmetric matrices}
Considering the sub-matrix in the form (\ref{equ-A_{i,f} => 0}) where $A_{n,i} \ne A_{i,n}$,
\begin{equation}
    A_{\mathrm{sub}} = \begin{bmatrix}
        A_{i,i} & A_{i,n} \\
        A_{n,i} & A_{n,n}
    \end{bmatrix},
\end{equation}
the current approach to factorize $A_{\mathrm{sub}}$ is to obtain the sparsification transformations $U_t$ and $L_t$ separately, where
\begin{equation}
\begin{aligned}
    U_t &= \begin{bmatrix}
        I & -T_2^\top \\
          & I
    \end{bmatrix}
    \begin{bmatrix}
        \Pi_{r_2}^\top \\ \Pi_{s_2}^\top
    \end{bmatrix},\\
    Lt &= \begin{bmatrix}
        \Pi_{r_1} & \Pi_{s_1}
    \end{bmatrix}
    \begin{bmatrix}
        I & \\
        -T_1 & I
    \end{bmatrix}.
\end{aligned}
\end{equation}
However, using unequal index sets ($s_1 \ne s_2$) often results in a poor condition number for $A_{i_{s_2},i_{s_1}}$(see \equ(\ref{equ-UtALt})) especially for a large depth level $L$, as it breaks the symmetry in the matrix-graph structure. A natural improvement, which is used in this paper, is to apply interpolative decomposition to $\begin{bmatrix}
    A_{n,i} \\ A_{i,n}^\top
\end{bmatrix}$, ensuring that $L_t = U_t^\top$. Specifically, we have 
\begin{equation}
    \begin{bmatrix}
        A_{n,i} \Pi_{r_1} \\ A_{i,n}^\top \Pi_{r_1}
    \end{bmatrix} = 
    \begin{bmatrix}
        A_{n,i} \Pi_{s_1} \\ A_{i,n}^\top \Pi_{s_1}
    \end{bmatrix} T_1 + \O(\eps),
\label{equ-ID [A_ni; A_in]}
\end{equation}
which turns out to be much more stable in the computation.

\section{Complexity Analysis}\label{comp_analy}
\subsection{Complexity of constructing tree structure}
Consider the tree structure with the finest level $L= \O(\log(N))$ (see \fig\ref{fig-split tree}), where the total number of vertices is $N$. We assume that the $i$th subgraph at level $l$, denoted by $G_{l,i}$, satisfies following properties:
\begin{enumerate}
\item Each subgraph $G_{l,i}$ consists of $N_l$ vertices with $N_{l} = \O(2^{-l} N)$. In particular, the size of $G_{L,i}$ at the leaf level is $\O(1)$.

\item Each subgraph $G_{l,i}$ attaches to a separator $S_{l,i}$ of size $c_l = \O({ (N_l)^{1/2} })$.

\item Each vertex in $G_{l,i}$ only has $\O(1)$ edges due to the local connectivity property.
\end{enumerate}

These properties are easily satisfied by the nested structure of finite element graph. In \algo\ref{alg-nested 2D}, the predominant time cost is finding separators, where each step checks $\O(1)$ neighbors to determine whether the neighbor belongs to $G_{l,i}$, and the expansion takes $c_l$ steps. Hence, the complexity for finding separators via binary search is given by:
\begin{equation}
	T_{l}^{\mathrm{Nd}} = \O({ (N_l)^{1/2} \log(N_l) }).
\end{equation}
Therefore, the overall  cost for constructing the tree structure of nested dissection is:
\begin{equation}
    T^{\mathrm{Nd}} = \sum_{l=1}^{L} 2^l T_{l}^{\mathrm{Nd}} 
    = \O({ \sum_{l=0}^{L-1} 2^l (N_l)^{1/2} \log(N_l) })
	= \O(N).
\end{equation}
\subsection{Complexity of factorization based on low rank approximations}\label{comp_analy2}
To estimate the cost of spaLU during factorization, we need the following sparsity assumption on the low rank structures:

\textit{Sparsity assumption}: let $e_l$ be the maximum size of segments before sparsification at level $l$. The size after sparsification, denoted by $e_l'$, satisfies the following property:
\begin{equation}\label{spaasump}
    e_l' \le q \cdot e_l + \sigma_e,
\end{equation}
where $q < 2^{-1/3}$ is the compression rate and $\sigma_e =\O(1)$ is a positive constant. 

\begin{remark}
Since the compression rate $q$ must be less than 1, our assumption that $q < 2^{-1/3}\approx 0.8$ is quite mild. It is verified by the numerical experiments in Section 6 for various PDEs, which implies the assumption is also practically reasonable.
\end{remark}

In the elimination of interiors, there are $2^L$ leaf subgraphs of size $N_L=\O(1)$, with neighboring sizes of $\O(1)$. Thus, the complexity of eliminating interiors is given by
\begin{eqnarray}
    T^{\mathrm{F,int}} = \O(2^L N_L^3) = \O(N).
\end{eqnarray}

As shown in \fig\ref{fig-segment_generate}, separator $S_{l,i}$ is split into at most $n_{j,l} = \O(2^{\frac{j-l}{2}})$ segments  $S_{l,i}^{j,1}$, $S_{l,i}^{j,2}$, $\dots$, at level $j$. This implies that segment $S_{l,i}^{j,k}$ has size at most
\begin{equation}
    s_j = \O({ \max_{0 \le l \le j} \frac{c_l}{n_{j,l}} })
    = \O({ 2^{-j/2} N^{1/2} }).
\end{equation}
Consider that there are $2^l$ separators $S_{l,1}, S_{l,2}, \dots$ at level $l$.
The total number of segments $S^{j,k}_{l,i}$ at level $j$, with $1 \le l \le j$, $1 \le i \le 2^l$, and $1 \le k \le n_{j,l}$, is
\begin{equation}
    n_j = \O( \sum_{l=1}^{j} 2^l \cdot n_{j,l} )
    = \O(2^j).
\end{equation}
Note that every segment is divided into finer regular segments every two levels, as shown in \fig\ref{fig-segment_generate}. Conversely, every segment $S^{j,k}_{l,i}$ at level $j$ will be merged after two levels of sparsification. In particular, considering $e_j$ at level $j$, since some segments are merged at level $j-1$, it holds $e_{j-1} \le 2 e_j' \le 2q e_j + 2 \sigma_e$.
At level $j-2$, a segment $S_{l,i}^{j-2,k}$ can be merged by some segments $S_{l,i}^{j-1,p}$ that were not merged at level $j-1$, i.e., there exists $h$ such that $S_{l,i}^{j-1,p} = S_{l,i}^{j,h}$, yielding $e_{j-2} \le 2q e_j + 2 \sigma_e$.
Given $e_L = N_L \in \O(1)$, we have:
\begin{equation}
    e_{L-2i+1} \le e_{L-2i} \le (2q)^{i} e_L + 2\sigma_e \sum_{j=0}^{i-1} (2q)^j = \O({(2q)^i})
\end{equation}
for $2 \le 2i \le L-1$, which implies $e_{L-i} = \O( {(2q)^{i/2}} )$ for $1 \le i < L$. Since the QR factorization dominates the computational time during the sparsification, the flops for sparsification at level $j$ is proportional to
\begin{equation}
    T^{\mathrm{F,sp}}_j = \O( e_j^2 e_j') = \O(e_j^3).
\end{equation}
Therefore, the complexity of total sparsification across all levels is bounded by
\begin{equation}
    T^{\mathrm{F,sp}} = \sum_{i=0}^{L-1} T^{\mathrm{F,sp}}_{L-i}
    = \O({ \sum_{i=0}^{L-1} 2^{L-i} (2q)^{3i/2} })
    = \O(N).      
\end{equation}
The time complexity of elimination, denoted by $T^{\mathrm{F,el}}$ can be derived in a manner analogous to $T^{\mathrm{F,sp}}$. Hence, the overall cost of spaLU is given by
\begin{equation}
    T^{\mathrm{F}}  = T^{\mathrm{F,int}} + T^{\mathrm{F,sp}} + T^{\mathrm{F,el}} = \O(N).
\end{equation}
Furthermore, since solving a triangular matrix equation of size $n \times n$ costs $\O(n^2)$, the time complexity of applying the spaLU factorization as a fast direct solver, derived in a similar way as before, is given by
\begin{equation}
    T^{\mathrm{S}} = \O({ \sum_{i=0}^{L-1} 2^{L-i} (2q)^i })
    = \O(N).
\end{equation}


\section{Numerical Experiments}
In this section, we test the performance of the algorithm by applying it to various PDE problems. All the equations are discretized by FEM using triangular meshes. We test the computational time and residuals with different degrees of freedom (DOF) $N$ for $\eps = 10^{-12}, 10^{-10}, 10^{-8}$, where $\eps$ is the tolerance defined in Section \ref{section-intro ID}. The algorithm is implemented in Matlab and all tests are conducted on a server with 128 GB RAM and an Intel Xeon CPU.

To illustrate the results, we use the following notations in the tables:
\begin{enumerate}
    \item $T^{\mathrm{Nd}}$ represents the time for constructing nested structure (in seconds).
    \item $T^{\mathrm{F}}$ refers to the factorization time via spaLU (in seconds).
    \item $T^{\mathrm{S}}$ is the time to apply the factorization to solving the linear system (in seconds).
    \item $R_{\mathrm{es}}$ denotes the relative residual $\|A x - b\|_2 / \|b\|_2$ after solving equation $\mathcal{A} x = b$.
    \item $\tilde{q}$ represents the numerical compression factor for segments sparsification, i.e., $\tilde{q} = \max_{l=1}^{L} \frac{e_l'}{e_l}$, which is used to  numerically verify the sparsity assumption $q$ in equation \eqref{spaasump}.
\end{enumerate}

\subsection{High-contrast Laplace equation}
\begin{figure}[htb]
    \centering
    \begin{subfigure}{0.48\textwidth}
        \centering
        \includegraphics [width=5.8cm] {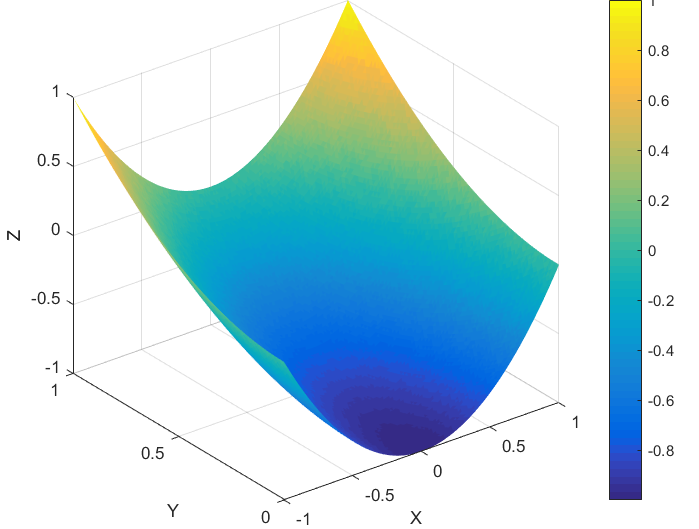}
        \caption{$\rho = 1$}
    \end{subfigure}
	\begin{subfigure}{0.48\textwidth}
		\centering
		\includegraphics [width=5.8cm] {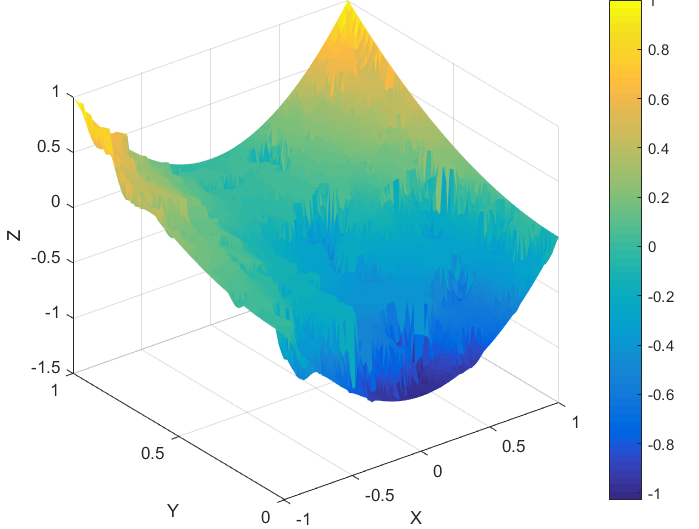}
		\caption{$\rho = 100$}
	\end{subfigure}
    \caption{Illustrations of solutions for Laplace equations in different media. (a) Constant medium. (b) High-contrast medium. 
    }
    \label{fig-experiment high-contrast a field}
\end{figure}

Consider the 2D Laplace equation:
\begin{equation}\label{high_lap}
\begin{cases}
    \nabla \cdot (a \nabla u) = f, & \Omega = [-1,1] \times [0,1], \\
    u|_{\partial \Omega} = g(x),
\end{cases}
\end{equation}
 where $f(x) = -4$, $g(x)  = x_1^2 + x_2^2 - 1$, and the variable coefficient $a(x)$ represents a high-contrast medium characterized by large values of $\rho$ and small values of $\rho^{-1}$, with $\rho\ge 1$. We generate $a$ by smoothing a random field in $[0,1]$, and take $a=\rho$ when the value of random field is greater than $0.5$ and  $a=\rho^{-1}$ otherwise. For $\rho = 1$, $a$ is simply a constant,  while $\rho = 100$, $a$ becomes highly irregular. Illustrations of the  solutions with $\rho = 1$ and $\rho = 100$ are provided in Figure \ref{fig-experiment high-contrast a field}.
One can observe there are significant roughness in the high-contrast medium.

\begin{figure}[htb]
    \centering
    \begin{subfigure}{0.49\textwidth}
        \centering
        \includegraphics [width=0.9\textwidth] {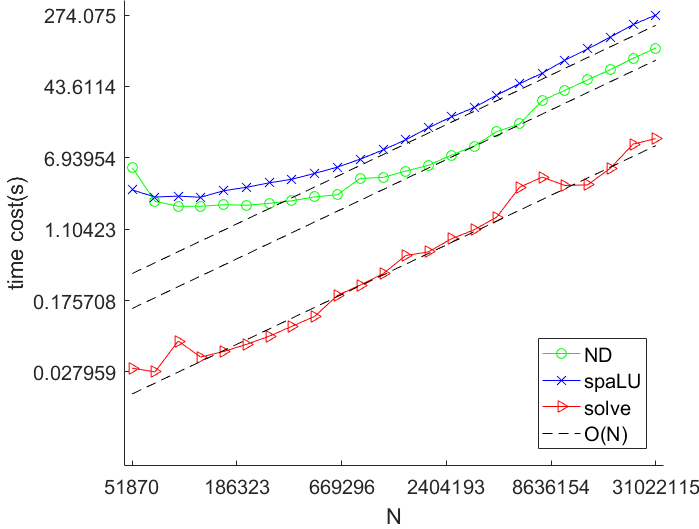}
        \caption{$\rho=1$}
    \end{subfigure}
    \begin{subfigure}{0.49\textwidth}
        \centering
        \includegraphics [width=0.9\textwidth] {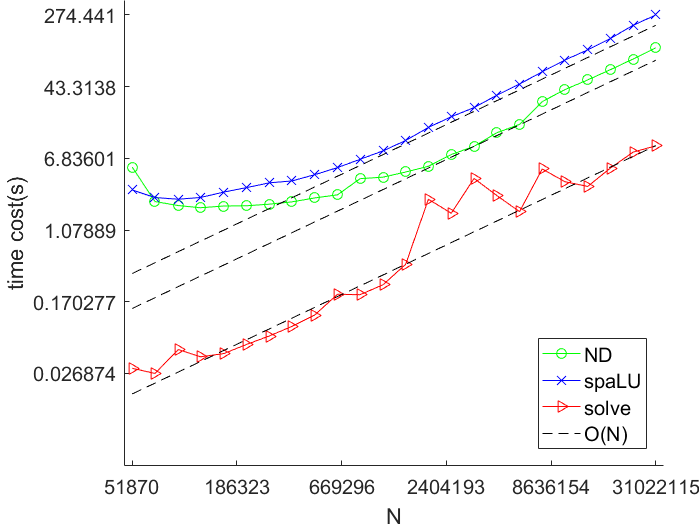}
        \caption{$\rho=100$}
    \end{subfigure}
    \caption{Solving the Laplace equation with FEM using spaLU:  the computational costs of constructing nested dissection data structure (ND), spaLU, and solving the linear system (solve) align closely with the expected linear complexity $\O(N)$, both in the constant medium case (a) and high-contrast medium case (b).}
    \label{fig-experiment high-contrast}
\end{figure}

Results for solving the discretized version of equation \eqref{high_lap} with $\eps = 10^{-12}$ are shown in Figure \ref{fig-experiment high-contrast}. Specifically, Figure \ref{fig-experiment high-contrast}(a) displays the results for the constant coefficient case $a = 1$. The three colored lines correspond to the computational time for constructing the nested dissection structure ($T^{\mathrm{Nd}}$), performing the recursive sparse LU decomposition ($T^{\mathrm{F}}$), and solving the linear system ($T^{\mathrm{S}}$), respectively. Similar results are shown in Figure \ref{fig-experiment high-contrast}(b) for $\rho = 100$. Both results demonstrate that the computational complexity is $\mathcal{O}(N)$, which agrees well with our theoretical analysis.

Detailed computational time and accuracy for various value of $N$ are presented in Table \ref{table1}. It can be observed that the performance of the solver is not significantly affected by the roughness of $a$. In particular, high-contrast medium slightly increases the computational cost and leads to worse residuals due to the larger condition number. One can also see that the numerical sparsity rate $\tilde{q}$ is only slightly larger for $\rho=100$ compared to $\rho=1$, which is expected since the rank of interpolative decomposition increases in high-contrast medium. However, $\tilde{q}$ is less than $0.8$ in all tests, supporting the reasonableness of our sparsity assumption in subsection \ref{comp_analy2}. It is also evident that the majority of the computational time is spent on constructing the tree structure for nested dissection ($T^{\mathrm{Nd}}$) and the recursive LU factorization ($T^{\mathrm{F}}$). Once the matrix has been factorized, the amount of time $T^{\mathrm{S}}$ used to solve the linear system is negligible. This is a key advantage of fast direct solver compared to iterative methods.


\begin{table}[htb]
    \begin{tabular}{lc|c|lll|c|c}
    \hline
        $N$ & $n \approx N^{1/2}$ & $\rho$ & $T^{\mathrm{Nd}}$ & $T^{\mathrm{F}}$ & $T^{\mathrm{S}}$ & $\tilde{q}$ & $R_{es}$ \\ \hline
        
        \multirow{2}{*}{364 514} & \multirow{2}{*}{604} & 1 & 2.29 & 3.99 & 0.09 & 0.63 & $6.58 \times 10^{-13}$ \\
        & & 100 & 2.23 & 3.82 & 0.090 & 0.75 & $8.74 \times 10^{-12}$ \\ \hline
        
        \multirow{2}{*}{1 110 389} & \multirow{2}{*}{1054} & 1 & 4.18 & 8.50 & 0.35 & 0.66 & $3.83 \times 10^{-12}$ \\
        & & 100 & 4.21 & 8.30 & 0.26 & 0.69 & $1.62 \times 10^{-11}$ \\ \hline
        
        \multirow{2}{*}{3 369 978} & \multirow{2}{*}{1836} & 1 & 9.37 & 25.32 & 1.08 & 0.67 & $1.09 \times 10^{-11}$ \\
        & & 100 & 9.31 & 25.26 & 4.07 & 0.70 & $2.05 \times 10^{-11}$ \\ \hline
        
        \multirow{2}{*}{10 225 605} & \multirow{2}{*}{3198} & 1 & 39.41 & 85.78 & 3.40 & 0.65 & $9.25 \times 10^{-12}$ \\
        & & 100 & 40.18 & 85.97 & 3.75 & 0.79 & $2.98 \times 10^{-11}$ \\ \hline
        
        \multirow{2}{*}{31 022 115} & \multirow{2}{*}{5570} & 1 & 117.11 & 274.07 & 11.37 & 0.66 & $1.36 \times 10^{-11}$ \\
        & & 100 & 117.64 & 274.44 & 9.44 & 0.77 & $3.40 \times 10^{-11}$\\ \hline
        
    \end{tabular}\caption{The computational time, sparsity rate and residuals for solving 2D Laplace equation using spaLU.}\label{table1}
\end{table}    

\subsection{Helmholtz equation in a regular domain}
Consider the 2D Helmholtz equation
\begin{align}
\begin{cases}
    \Delta u+ k^2 u = f, \quad \Omega = [-1,1] \times [0,1], \\
    u|_{\partial \Omega} = g(x),
    \end{cases}
\end{align}
with $k = \sqrt{2}$,  $f(x) = -1$ and $g(x)  = e^{x_1+x_2}$. The equation is discretized by FEM with linear elements on a triangular mesh. 

\begin{figure}[htb]
    \centering
    \begin{subfigure}{0.48\textwidth}
        \centering
        \includegraphics [width=0.9\textwidth] {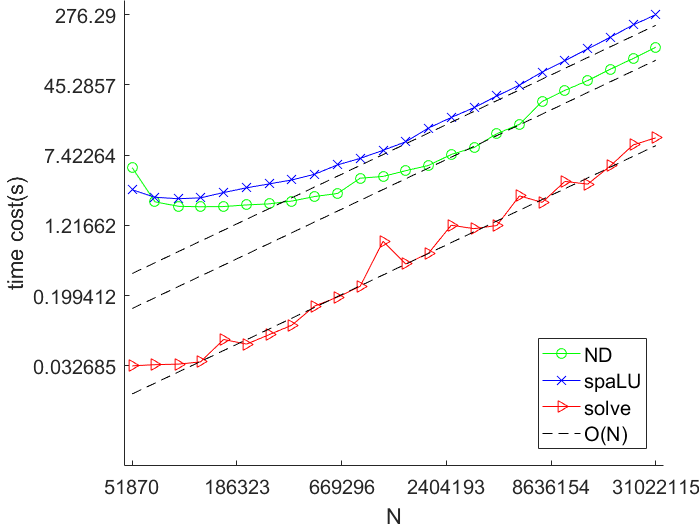}
        \caption{$\eps=10^{-12}$}
    \end{subfigure}
    \begin{subfigure}{0.48\textwidth}
        \centering
        \includegraphics [width=0.9\textwidth] {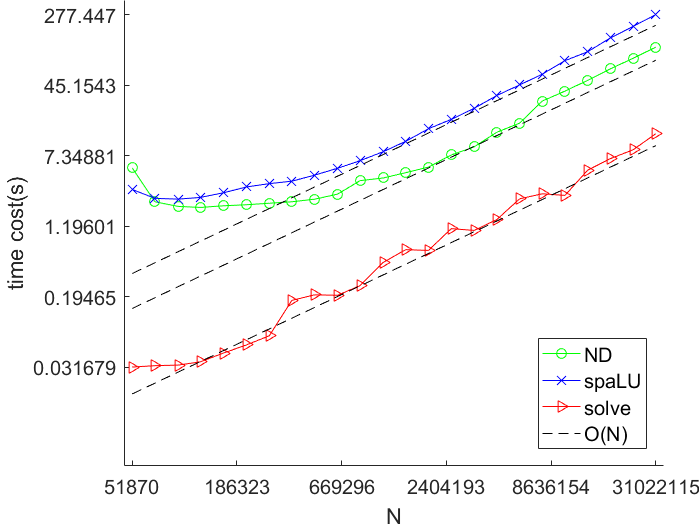}
        \caption{$\eps=10^{-10}$}
    \end{subfigure}
    \\
    \begin{subfigure}{0.48\textwidth}
        \centering
        \includegraphics [width=0.9\textwidth] {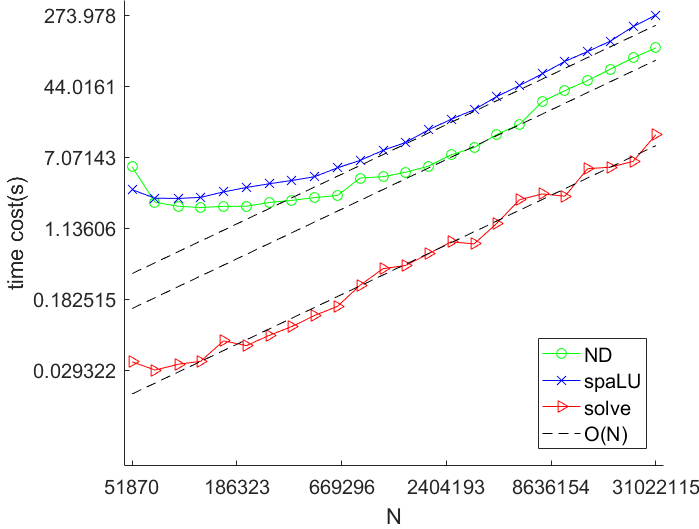}
        \caption{$\eps=10^{-8}$}
    \end{subfigure}
   \caption{Solving the discretized Helmholtz equation with  up to 31 million unknowns via spaLU. The computational time for nested dissection tree (ND), factorization (spaLU), and solving the linear systems based on the factorization (solve) are plotted. The black dashed line indicates the computational complexity of all the three processes is on the order of $\mathcal{O}(N)$.}
    \label{fig-experiment helmholtz}
\end{figure}

Figure \ref{fig-experiment helmholtz} shows the computational results for solving the linear system with DOF up to 31 million. It is clear that the computational complexity is on the order of $\mathcal{O}(N)$ for all three processes, namely, tree structure construction, factorization and solving, which is consistent with our theoretical analysis. The figure also shows that reducing the tolerance only slightly decreases the computational cost. The reason is due to the rapid decay of singular values of the submatrices associated with segments, leading to a relatively stable size of skeletons even with smaller tolerance.  More detailed results are provided in Table \ref{table2}. One can see the sparsity rate $\tilde{q}$ is almost constant (around $0.6$) across different DOF and tolerance levels, implying that it may be an intrinsic property related to the PDE itself. The residual $R_{es}$ closely matches the compression tolerance $\eps$, with an average loss of one digit, indicating the algorithm is very stable. Furthermore, the time required to apply the factorization for solving the linear system is significantly less than the time needed for the factorization process itself, as observed in the previous example. This efficiency suggests that our solver is highly effective when multiple right hand sides need to be solved.

\begin{table}
    \begin{tabular}{lc|c|lll|c|c}
    \hline
        $N$ & $n \approx N^{1/2}$ & $\eps$ & $T^{\mathrm{Nd}}$ & $T^{\mathrm{F}}$ & $T^{\mathrm{S}}$ & $\tilde{q}$ & $R_{es}$ \\ \hline
        
        \multirow{3}{*}{364 514} & \multirow{3}{*}{604} & $10^{-8}$ & 2.31 & 3.87 & 0.09 & 0.66 & $2.50 \times 10^{-7}$ \\
        & & $10^{-10}$ & 2.23 & 3.79 & 0.17 & 0.63 & $7.28 \times 10^{-10}$ \\
        & & $10^{-12}$ & 2.28 & 3.92 & 0.093 & 0.63 & $7.01 \times 10^{-13}$ \\ \hline
        
        \multirow{3}{*}{1 110 389} & \multirow{3}{*}{1054} & $10^{-8}$ & 4.26 & 8.40 & 0.40 & 0.67 & $5.61 \times 10^{-7}$ \\
        & & $10^{-10}$ & 4.12 & 8.10 & 0.47 & 0.68 & $1.16 \times 10^{-9}$ \\
        & & $10^{-12}$ & 4.29 & 8.41 & 0.81 & 0.66 & $3.95 \times 10^{-12}$ \\ \hline
        
        \multirow{3}{*}{3 369 978} & \multirow{3}{*}{1836} & $10^{-8}$ & 9.04 & 23.88 & 0.75 & 0.68 & $8.11 \times 10^{-7}$ \\
        & & $10^{-10}$ & 9.17 & 24.72 & 1.06 & 0.67 & $1.46 \times 10^{-9}$ \\
        & & $10^{-12}$ & 9.03 & 25.24 & 1.12 & 0.67 & $1.07 \times 10^{-11}$ \\ \hline
        
        \multirow{3}{*}{10 225 605} & \multirow{3}{*}{3198} & $10^{-8}$ & 39.20 & 83.14 & 2.58 & 0.65 & $1.16 \times 10^{-6}$ \\
        & & $10^{-10}$ & 38.64 & 84.52 & 2.63 & 0.65 & $2.12 \times 10^{-9}$ \\
        & & $10^{-12}$ & 39.36 & 84.66 & 3.77 & 0.65 & $9.23 \times 10^{-12}$ \\ \hline
        
        \multirow{3}{*}{31 022 115} & \multirow{3}{*}{5570} & $10^{-8}$ & 117.38 & 273.98 & 12.56 & 0.66 & $1.49 \times 10^{-6}$ \\
        & & $10^{-10}$ & 119.84 & 277.44 & 12.89 & 0.66 & $2.68 \times 10^{-9}$ \\
        & & $10^{-12}$ & 119.86 & 276.28 & 11.60 & 0.66 & $1.35 \times 10^{-11}$ \\
        \hline
    \end{tabular}
    \caption{The computational time, sparsity rate and residuals for solving Helmholtz equation in a regular domain using spaLU.}\label{table2}
\end{table}    

\subsection{Helmholtz equation in an irregular domain}
\begin{figure}[htb]
    \centering
    \begin{subfigure}{0.48\textwidth}
        \centering
            \includegraphics[width=0.9\textwidth]{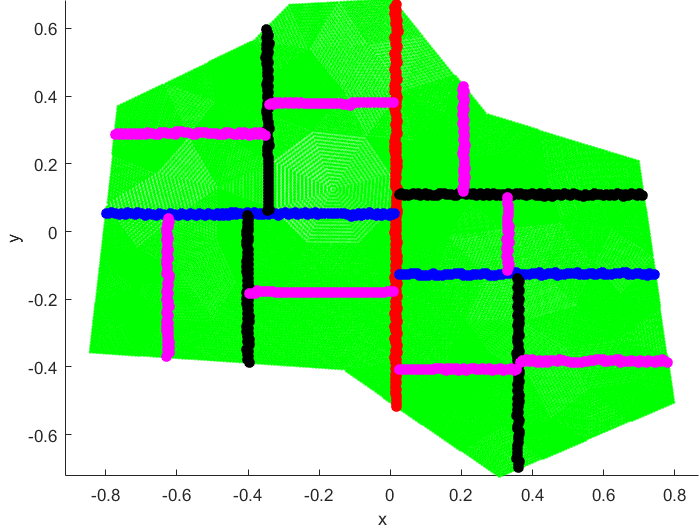}
        \caption{Nested partition for an irregular domain}
    \end{subfigure}
    \begin{subfigure}{0.48\textwidth}
        \centering
        \includegraphics[width=0.9\textwidth] {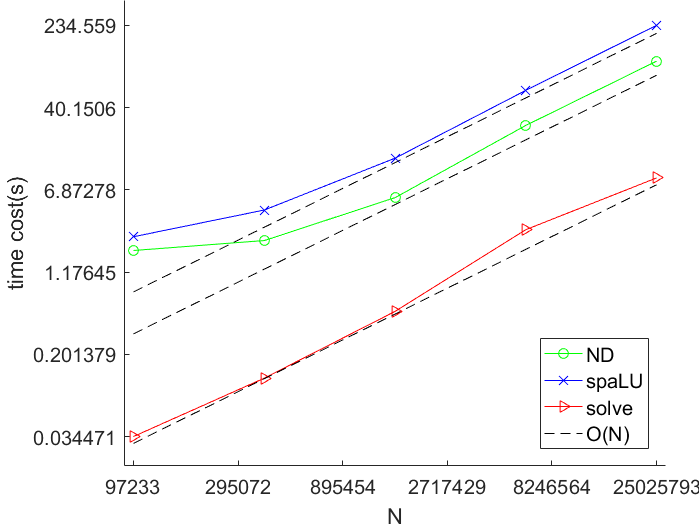}
        \caption{Computational complexity versus the DOF}
    \end{subfigure}
    \caption{Solving the discretized Helmholtz equation in an irregular domain via spaLU. The computational time for nested dissection tree (ND), factorization (spaLU), and solving the linear systems based on the factorization (solve) are plotted. The black dashed line indicates the computational complexity of all the three processes is on the order of $\mathcal{O}(N)$.}
    \label{fig-experiment irregular}
\end{figure}

Consider the 2D Helmholtz equation within an irregular domain $\Omega$:
\begin{equation}
\begin{cases}
    \Delta u + k^2 u = f(x), \quad 
    \mbox{for } x \in \Omega, \\
    u|_{\partial \Omega} = g(x),
    \end{cases}
\end{equation}
where $k = \sqrt{2}$, $f(x) = -1$, $g(x)  = e^{x_1+x_2}$, and $\Omega$ is a polygonal domain, as shown in \fig\ref{fig-experiment irregular}(a). The figure also shows the structure of nested partition. Computational results for $\eps=10^{-12}$ are shown in \fig\ref{fig-experiment irregular}(b). Once again, we can see that all the three processes, including the tree structure construction, recursive sparse LU decomposition and applying the factorization for solving linear system, are all on the order of $\O(N)$. Compared to the regular domain, the compression rate $\tilde{q}$ slightly increases, but is still less than $0.8$ and satisfies our sparsity assumption. In terms of efficiency, the factorization time in the irregular domain is roughly $0.106$ million DOF per second, compared to approximately $0.112$ million DOF per second in the regular domain. The difference is minor, which implies the solver is efficient in the irregular domain as well.

\begin{table}
    \begin{tabular}{lc|lll|c|c}
     \hline
        $N$ & $n \approx N^{1/2}$ & $T^{\mathrm{Nd}}$ & $T^{\mathrm{F}}$ & $T^{\mathrm{S}}$ & $\tilde{q}$ & $R_{es}$ \\ \hline
        
        97233 & 312 & 1.86 & 2.50 & 0.03 & 0.72 & $6.32 \times 10^{-13}$ \\ \hline
        390049 & 625 & 2.31 & 4.45 & 0.12 & 0.55 & $8.37 \times 10^{-13}$ \\ \hline
        1562433 & 1250 & 5.81 & 13.55 & 0.51 & 0.76 & $9.20 \times 10^{-13}$ \\ \hline
        6254209 & 2501 & 27.36 & 58.06 & 2.92 & 0.78 & $1.43 \times 10^{-12}$ \\ \hline
        25025793 & 5003 & 108.53 & 234.55 & 8.86 & 0.79 & $1.82 \times 10^{-12}$ \\ \hline    
    \end{tabular}\caption{The computational time, sparsity rate and residuals for solving Helmholtz equation in a irregular domain using spaLU.}\label{table3}
\end{table}

\subsection{Laplace equation with anisotropic coefficient}
Consider the 2D Laplace equation in anisotropic medium
\begin{equation}\label{anis_lap}
\begin{cases}
    \nabla \cdot (D \nabla u(x)) = f, \quad \Omega = [-1,1] \times [0,1], \\
    u|_{\Gamma_1} = g(x), \\
    \frac{\partial u}{\partial n}|_{\Gamma_2} = h(x),
    \end{cases}
\end{equation}
where  $f(x) = -4$, $g(x)  = x_1^2 + x_2^2 - 1$ and $h(x)  = 2x_2$, with $\Gamma_2 = \left\lbrace (x_1, x_2) \mid -1 < x_1 < 1, \; x_2 = 1 \right\rbrace$ and $\Gamma_1 = \partial \Omega \setminus \Gamma_2$. The coefficient $D=\begin{bmatrix}1 & 1 \\ 0 & 1\end{bmatrix}$ is an asymmetric matrix, representing the anisotropic medium. The resulting stiffness matrix $\mathcal{A}$ from FEM is unsymmetric in this case. 

\begin{figure}[htb]
	\centering
	\begin{subfigure}{0.48\textwidth}
		\centering
		\includegraphics [width=0.9\textwidth] {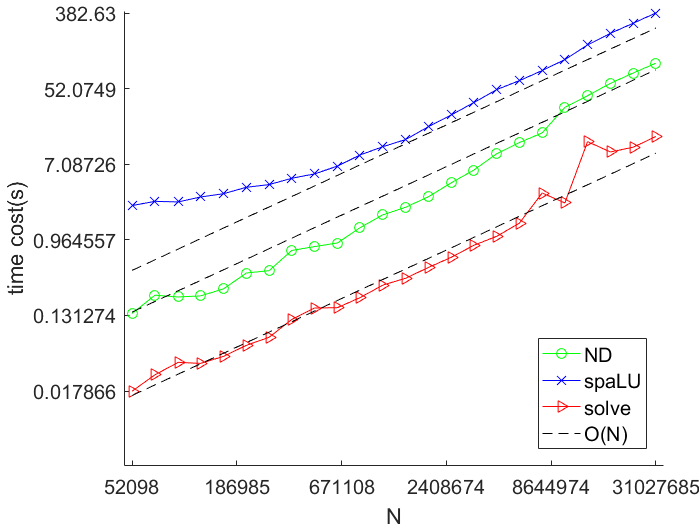}
		\caption{$\eps=10^{-12}$}
	\end{subfigure}
	\begin{subfigure}{0.48\textwidth}
		\centering
		\includegraphics [width=0.9\textwidth] {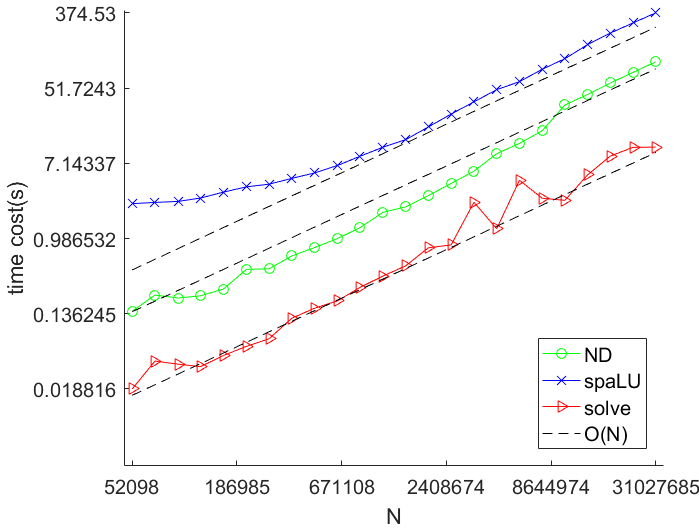}
		\caption{$\eps=10^{-10}$}
	\end{subfigure}
	\\
	\begin{subfigure}{0.48\textwidth}
		\centering
		\includegraphics [width=0.9\textwidth] {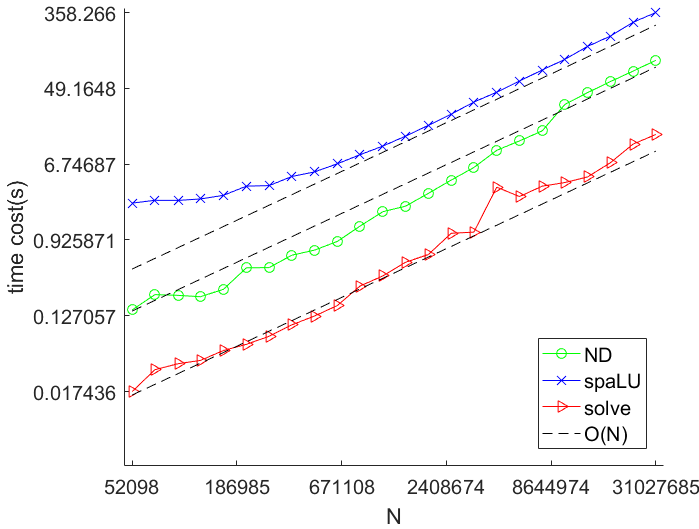}
		\caption{$\eps=10^{-8}$}
	\end{subfigure}
	\caption{Solving the discretized anisotropic Laplace equation via spaLU. The computational time for constructing nested dissection tree (ND), factorization (spaLU), and solving the linear systems based on the factorization (solve) are plotted. The black dashed line indicates the computational complexity of all the three processes is on the order of $\mathcal{O}(N)$.}
	\label{fig-experiment anisotropy}
\end{figure}

Numerical results for numerically solving equation \eqref{anis_lap} with up to 31 million DOF via spaLU are shown in \fig\ref{fig-experiment anisotropy}. Considering that Cholesky factorization for symmetric matrices generally requires about half the complexity of LU factorization for unsymmetric matrices of the same size,  the factorization for an unsymmetric matrix typically takes twice as long as compared to the same sized symmetric matrix. However, our numerical experiments show that for unsymmetric matrices, spaLU is significantly faster than twice the time required of symmetric matrices. In particular, compared to the computational time for isotropic case in Example 1, where symmetric matrices are factorized, the factorization time for the anisotropic case only increases by roughly $30$ percent for the same amount of DOF. It implies the proposed sampling method considerably accelerate the interpolative decomposition for unsymmetric matrices. Meanwhile, the method still follows the rule of linear complexity in the unsymmetric case, with the sparsity rate $\tilde{q}$ is consistently less than $0.8$, thus satisfying our sparsity assumption.

\begin{table}
	\begin{tabular}{lc|c|lll|c|c}
     \hline
		$N$ & $n \approx N^{1/2}$ & $\eps$ & $T^{\mathrm{Nd}}$ & $T^{\mathrm{F}}$ & $T^{\mathrm{S}}$ & $\tilde{q}$ & $R_{es}$ \\ \hline
		
		\multirow{3}{*}{365 118} & \multirow{3}{*}{604} & $10^{-8}$ & 0.62 & 4.88 & 0.10 & 0.73 & $5.96 \times 10^{-7}$ \\
		& & $10^{-10}$ & 0.62 & 4.80 & 0.12 & 0.57 & $2.02 \times 10^{-10}$ \\
		& & $10^{-12}$ & 0.73 & 4.89 & 0.12 & 0.69 & $1.02 \times 10^{-12}$ \\ \hline
		
		\multirow{3}{*}{1 111 443} & \multirow{3}{*}{1054} & $10^{-8}$ & 1.94 & 10.85 & 0.36 & 0.76 & $1.16 \times 10^{-6}$ \\
		& & $10^{-10}$ & 1.96 & 10.91 & 0.37 & 0.60 & $2.32 \times 10^{-10}$ \\
		& & $10^{-12}$ & 1.87 & 11.28 & 0.29 & 0.74 & $2.80 \times 10^{-12}$ \\ \hline
		
		\multirow{3}{*}{3 371 814} & \multirow{3}{*}{1836} & $10^{-8}$ & 6.19 & 34.33 & 1.14 & 0.73 & $1.54 \times 10^{-6}$ \\
		& & $10^{-10}$ & 5.81 & 35.96 & 2.56 & 0.62 & $3.76 \times 10^{-10}$ \\
		& & $10^{-12}$ & 6.09 & 36.08 & 0.84 & 0.73 & $4.10 \times 10^{-11}$ \\ \hline
		
		\multirow{3}{*}{10 228 803} & \multirow{3}{*}{3198} & $10^{-8}$ & 32.24 & 105.15 & 4.16 & 0.69 & $1.01 \times 10^{-6}$ \\
		& & $10^{-10}$ & 33.13 & 111.92 & 2.70 & 0.61 & $4.57 \times 10^{-10}$ \\
		& & $10^{-12}$ & 31.74 & 111.85 & 2.56 & 0.68 & $2.89 \times 10^{-12}$ \\ \hline
		
		\multirow{3}{*}{31 027 685} & \multirow{3}{*}{5570} & $10^{-8}$ & 102.40 & 358.27 & 14.65 & 0.73 & $1.05 \times 10^{-6}$ \\
		& & $10^{-10}$ & 103.27 & 374.53 & 10.94 & 0.62 & $6.28 \times 10^{-10}$ \\
		& & $10^{-12}$ & 101.79 & 382.63 & 14.72 & 0.69 & $2.94 \times 10^{-12}$ \\ \hline
		
	\end{tabular}\caption{The computational time, sparsity rate and residuals for solving anisotropic Laplace equation using spaLU.}\label{table4}
\end{table}
	
\section{Conclusion}
In this paper, we have demonstrated an efficient approach for solving large sparse linear systems resulting from the FD or FEM discretization of PDEs, based on the proposed recursive sparse LU decomposition.
By employing nested dissection and low rank approximations, we reorganized the matrix structure and successfully implemented a hierarchical sparsification of dense blocks on the separators.
Our hybrid sampling algorithm, which combines ideas from randomized sampling and the FMM, effectively skeletonized the dense blocks, enabling the construction of a fast direct solver applicable to both  symmetric and unsymmetric matrices.
The theoretical complexity of the solver is $\O(N)$ under a mild compression rate assumption, which has been confirmed by the numerical experiments.
Potential applications of the solver include large-scale optimization problems with PDE constraint or inverse problems, where a large number of forward solves are required.
Extension of the algorithm to three dimensional PDEs is currently under investigation.

\bibliographystyle{plain} 
\bibliography{references3}

\begin{thebibliography}{10}

\bibitem{HODLR}
Amirhossein Aminfar, Sivaram Ambikasaran, and Eric Darve.
\newblock A fast block low-rank dense solver with applications to
  finite-element matrices.
\newblock {\em Journal of Computational Physics}, 304:170--188, 2016.

\bibitem{Intro-solveEq}
Mario Bebendorf.
\newblock Efficient inversion of the galerkin matrix of general second-order
  elliptic operators with nonsmooth coefficients.
\newblock {\em Mathematics of computation}, 74(251):1179--1199, 2005.

\bibitem{spaND}
L{\'e}opold Cambier, Chao Chen, Erik Boman, Sivasankaran Rajamanickam, Raymond
  Tuminaro, and Eric Darve.
\newblock An algebraic sparsified nested dissection algorithm using low-rank
  approximations.
\newblock {\em SIAM Journal on Matrix Analysis and Applications},
  41(2):715--746, 2020.

\bibitem{HSS}
Shiv Chandrasekaran, Ming Gu, and Timothy Pals.
\newblock A fast ${ULV}$ decomposition solver for hierarchically semiseparable
  representations.
\newblock {\em SIAM Journal on Matrix Analysis and Applications},
  28(3):603--622, 2006.

\bibitem{ID-LowRank}
Hongwei Cheng, Zydrunas Gimbutas, Per-Gunnar Martinsson, and Vladimir Rokhlin.
\newblock On the compression of low rank matrices.
\newblock {\em SIAM Journal on Scientific Computing}, 26(4):1389--1404, 2005.

\bibitem{ND-intro_lite}
Timothy Davis.
\newblock {\em Direct Methods for Sparse Linear Systems}.
\newblock SIAM, 2006.

\bibitem{ND-NPC}
Michael Garey and David Johnson.
\newblock {\em Computers and intractability: A guide to the theory of
  np-completeness}, volume~24.
\newblock Society for Industrial and Applied Mathematics, 1982.

\bibitem{ND-intro_first}
Alan George.
\newblock Nested dissection of a regular finite element mesh.
\newblock {\em SIAM journal on numerical analysis}, 10(2):345--363, 1973.

\bibitem{MinDegreeOrder}
Alan George and Joseph Liu.
\newblock The evolution of the minimum degree ordering algorithm.
\newblock {\em Siam review}, 31(1):1--19, 1989.

\bibitem{MatrixComputation}
Gene Golub and Charles Van~Loan.
\newblock {\em Matrix computations}.
\newblock Johns Hopkins University Press, 4th edition, 2013.

\bibitem{FMM-lecture}
Leslie Greengard and Vladimir Rokhlin.
\newblock A fast algorithm for particle simulations.
\newblock {\em Journal of computational physics}, 73(2):325--348, 1987.

\bibitem{Strong_RRQR}
Ming Gu and Stanley Eisenstat.
\newblock Efficient algorithms for computing a strong rank-revealing qr
  factorization.
\newblock {\em SIAM Journal on Scientific Computing}, 17(4):848--869, 1996.

\bibitem{Multigrid}
Wolfgang Hackbusch.
\newblock {\em Multi-grid methods and applications}.
\newblock Springer Berlin Heidelberg, Berlin, Heidelberg, 1985.

\bibitem{Hierarchical}
Wolfgang Hackbusch.
\newblock {\em Hierarchical Matrices: Algorithms and Analysis}, chapter
  Introduction, pages 3--24.
\newblock Springer Berlin Heidelberg, Berlin, Heidelberg, 2015.

\bibitem{H^2}
Wolfgang Hackbusch and Steffen B{\"o}rm.
\newblock ${H}^2$-matrix approximation of integral operators by interpolation.
\newblock {\em Applied numerical mathematics}, 43(1-2):129--143, 2002.

\bibitem{RandomMatrix}
Nathan Halko, Per-Gunnar Martinsson, and Joel Tropp.
\newblock Finding structure with randomness: Probabilistic algorithms for
  constructing approximate matrix decompositions.
\newblock {\em SIAM review}, 53(2):217--288, 2011.

\bibitem{cave2014}
Jun Lai, Sivaram Ambikasaran, and Leslie Greengard.
\newblock A fast direct solver for high frequency scattering from a large
  cavity in two dimensions.
\newblock {\em SIAM Journal on Scientific Computing}, 36(6):B887--B903, 2014.

\bibitem{Lai2015194}
Jun Lai, Motoki Kobayashi, and Alex Barnett.
\newblock A fast and robust solver for the scattering from a layered periodic
  structure containing multi-particle inclusions.
\newblock {\em Journal of Computational Physics}, 298:194 -- 208, 2015.

\bibitem{inverseP2022}
Jun Lai and Jinrui Zhang.
\newblock Fast inverse elastic scattering of multiple particles in three
  dimensions.
\newblock {\em Inverse Problems}, 38(10):104002, 2022.

\bibitem{ND-timecost}
Richard Lipton, Donald Rose, and Robert Tarjan.
\newblock Generalized nested dissection.
\newblock {\em SIAM journal on numerical analysis}, 16(2):346--358, 1979.

\bibitem{MaLin2024}
Haoran Ma, Gang Bao, Jun Lai, and Junshan Lin.
\newblock Inverse design of a grating metasurface for enhancing spontaneous
  emission through hyperbolic metamaterials.
\newblock {\em J. Opt. Soc. Am. B}, 41(2):A79--A85, Feb 2024.

\bibitem{FMM-skeleton}
Victor Minden, Kenneth Ho, Anil Damle, and Lexing Ying.
\newblock A recursive skeletonization factorization based on strong
  admissibility.
\newblock {\em Multiscale Modeling \& Simulation}, 15(2):768--796, 2017.

\bibitem{ND-graph_for_eliminate}
Seymour Parter.
\newblock The use of linear graphs in gauss elimination.
\newblock {\em SIAM review}, 3(2):119--130, 1961.

\bibitem{FMM-lecture2}
Vladimir Rokhlin.
\newblock Rapid solution of integral equations of scattering theory in two
  dimensions.
\newblock {\em Journal of Computational Physics}, 86(2):414--439, 1990.

\bibitem{GMRES}
Youcef Saad and Martin Schultz.
\newblock Gmres: A generalized minimal residual algorithm for solving
  nonsymmetric linear systems.
\newblock {\em SIAM Journal on Scientific and Statistical Computing},
  7(3):856--869, 1986.

\bibitem{ILU}
Yousef Saad.
\newblock ${ILUT}$: A dual threshold incomplete ${LU}$ factorization.
\newblock {\em Numerical linear algebra with applications}, 1(4):387--402,
  1994.

\end{thebibliography}

\end{document}